\documentclass[12pt,reqno]{amsart}

\usepackage{enumerate}
\usepackage{amsrefs}
\usepackage{amsfonts,amsmath,amssymb,amsthm,bm}
\usepackage{mathrsfs}
\usepackage[T1]{fontenc}
\usepackage{booktabs,array}
\usepackage[margin=1in]{geometry}
\usepackage{url}
\usepackage[
linktocpage=true,
colorlinks,
citecolor=magenta,
linkcolor=blue,
urlcolor=magenta
]{hyperref}

\usepackage{tikz}

\allowdisplaybreaks
\numberwithin{equation}{section}
\newtheorem{thm}{Theorem}[section]
\newtheorem{cor}[thm]{Corollary}
\newtheorem{lem}[thm]{Lemma}
\newtheorem{prop}[thm]{Proposition}

\theoremstyle{definition}

\newtheorem{prob}[thm]{Problem}

\theoremstyle{remark}
\newtheorem{rem}[thm]{Remark}

\newtheorem*{ack}{Acknowledgments}

\newcommand{\dd}{\,d}

\hypersetup{
 pdftitle={Brunn-Minkowski inequalities for chord integrals},
 pdfsubject={Chord integrals, convex lifting, covariograms and Riesz energies}
}
\begin{document}

\title[Brunn-Minkowski inequalities for chord integrals]
{Brunn-Minkowski inequalities for chord integrals}

\author{Haizhong Li}
\address{Department of Mathematical Sciences, Tsinghua University, Beijing 100084, P.R. China}
\email{\href{mailto:lihz@tsinghua.edu.cn}{lihz@tsinghua.edu.cn}}

\author{Yao Wan}
\address{Department of Mathematics, The Chinese University of Hong Kong, Shatin, N.T., Hong Kong}
\email{\href{mailto:yaowan@cuhk.edu.hk}{yaowan@cuhk.edu.hk}}

\author{Yijia Zhang}
\address{Department of Mathematical Sciences, Tsinghua University, Beijing 100084, P.R. China}
\email{\href{mailto:zhang-yj23@mails.tsinghua.edu.cn}{zhang-yj23@mails.tsinghua.edu.cn}}
%\date{}

\keywords{Chord integral; Brunn-Minkowski inequality; dual quermassintegral; chord Minkowski problem; affine-section integral; Riesz energy}

\subjclass[2020]{52A40; 52A22; 26D15}

\begin{abstract}
Chord integrals are moments of chord lengths and translation-invariant averages of dual quermassintegrals. In this paper, we prove their Brunn-Minkowski inequality for $q\in[1,n+2]\cup\{n+3\}$, with equality precisely for homothetic bodies. 
A transfer from symmetric dual quermassintegrals to covariogram level sets gives, in dimension two, the full range $1\leq q\leq5$. 
For every $0<q<1$, we establish the same inequality for simplices and an explicit factor version for all convex bodies.
On the other hand, rectangular perturbations give counterexamples for every $q>n+1+\beta_*$, where $\beta_*\in(2,3)$.
As consequences, we obtain Minkowski inequalities and uniqueness results for chord measures. We also prove Brunn-Minkowski inequalities for affine-section integrals and mixed Riesz energies.
\end{abstract}

\maketitle

\section{Introduction}

The Brunn-Minkowski theory studies geometric invariants of convex bodies through their behavior under Minkowski addition. Its classical invariants include the quermassintegrals and their associated area measures; see, for example, \cite{Schneider}.
The dual Brunn-Minkowski theory, initiated by Lutwak \cite{Lutwak75}, is built around the dual quermassintegrals
\[
 \widetilde V_r(K,z)
 =\frac{1}{n}\int_{\mathbb S^{n-1}}\rho_{K,z}(u)^r\,\dd u,
\]
where $\rho_{K,z}$ denotes the radial function of $K$ with respect to $z$. We denote $\widetilde V_r(K)=\widetilde V_r(K,0)$.
Dual quermassintegrals are generally origin-dependent, whereas the quermassintegrals are translation invariant. The study of these two families, together with their associated geometric measures, forms a central part of modern convex geometric analysis.

Integral geometry provides another natural source of global invariants by averaging over affine subspaces with respect to their invariant measures; see, e.g., \cites{Santalo,Ren,SchneiderWeil}.
For a convex body $K$, integrating over the lines that meet $K$ naturally produces the length of a random chord. Its power moments give rise to the \textit{chord integrals} studied in this paper. 
These functionals admit several equivalent interpretations: they can be written in terms of $X$-ray transforms, related through dual kinematic formulas to averages of affine sections, and represented by Riesz-type integrals of the characteristic function of $K$; see \cites{Gardner,Koldobsky,GZ98,Zhang99,LXYZ}.

Let $\mathcal K^n$ denote the family of convex bodies in $\mathbb R^n$, $n\geq 2$, and let $\mathcal A(n,1)$ be the affine Grassmannian of lines. 
We normalize its invariant measure so that it is a probability measure on directions and Lebesgue measure on the orthogonal displacements. For $q\geq 0$, the $q$-th chord integral of $K\in \mathcal K^n$ is
\begin{equation}\label{def-Iq}
 I_q(K)
 =\int_{\mathcal A(n,1)} |K\cap\ell|^q\,\dd\ell,
\end{equation}
where the integrand is understood as $\chi_{\{\ell\cap K\neq\varnothing\}}$ for $q=0$.
Thus, $I_q(K)$ is the unnormalized $q$-th moment of chord length with respect to the invariant measure on affine lines. The functional $I_q$ is invariant under translations and rotations and is homogeneous of degree $n+q-1$.
Although we define $I_q$ above for $q\geq0$, the same chord integral is finite for every $q>-1$; see \cite[Lemma~3.1]{LXYZ}.

To describe the chord integral more explicitly, let $\omega_j$ denote the volume of the $j$-dimensional Euclidean unit ball, let $\dd\sigma$ be the probability measure on $\mathbb S^{n-1}$, and define the $X$-ray function
\[
 X_K(y,u)=|K\cap(y+\mathbb Ru)|, \qquad y\in u^\perp,\quad u\in\mathbb S^{n-1}.
\]
Then
\begin{equation}\label{intro-Xray}
 I_q(K)
 = \int_{\mathbb S^{n-1}}\int_{K|u^\perp} X_K(y,u)^q\,\dd y\,\dd\sigma(u),\qquad q>-1,
\end{equation}
where $K|u^\perp$ is the image of the orthogonal projection of $K$ onto $u^\perp$; see \cite[(3.2)--(3.3)]{LXYZ}. This representation connects chord integrals with $X$-ray transforms and geometric tomography. 
In particular, $I_2(K)$ is the total squared $X$-ray signal of $K$.

There is also a natural connection with the dual Brunn-Minkowski theory. 
A chord calculation gives
\begin{equation}\label{intro-dual}
 I_q(K)
 = \frac{q}{\omega_n}\int_K \widetilde V_{q-1}(K,x)\,\dd x,\qquad q>0;
\end{equation}
see \cites{GZ98,Zhang99}.
Hence, chord integrals may be viewed as translation-invariant averages of the origin-dependent dual quermassintegrals. This identity is also the starting point for the construction of chord measures.

For $q>1$, the chord integral admits the following Riesz-energy and covariogram representations; see, e.g., \cites{Ren, Santalo, LXYZ, HYZ}.
For $K\in\mathcal{K}^n$, its covariogram function is given by
\begin{align*}
    g_K(z)=|K\cap(K+z)|=\chi_K*\chi_{-K}(z),
\end{align*}
where $\chi_K*\chi_{-K}$ denotes the convolution of the indicator functions $\chi_K$ and $\chi_{-K}$. Then
\begin{equation}\label{intro-cov-formula}
 \begin{split}
 I_q(K)
 &= \frac{q(q-1)}{n\omega_n} \int_K\int_K |x-y|^{q-n-1}\,\dd x\,\dd y\\
 &= \frac{q(q-1)}{n\omega_n} \int_{\mathbb R^n}g_K(z)|z|^{q-n-1}\,\dd z.
 \end{split}
\end{equation}
% For $1<q<n+1$, this is a Riesz self-energy, whereas for
% $q\geq n+1$ it can be viewed as a positive moment of pairwise distances. 
For $0<q<1$, the corresponding interior energy diverges, and the finite representation is instead given by the covariogram deficit; see \cites{Lud,Qin}:
\begin{equation}\label{intro-fractional}
 \begin{split}
 I_q(K)
 &=\frac{q(1-q)}{n\omega_n} \int_K\int_{K^c}|x-y|^{q-n-1}\,\dd y\,\dd x\\
 &=\frac{q(1-q)}{n\omega_n} \int_{\mathbb R^n} \bigl(|K|-g_K(z)\bigr)|z|^{q-n-1}\,\dd z.
 \end{split}
\end{equation}
Consequently, up to a constant, $I_q$ is the fractional perimeter of order $1-q$:
\[
 I_q(K)=\frac{q(1-q)}{n\omega_n}P_{1-q}(K),
 \qquad
 P_s(K) =\int_K\int_{K^c} \frac{\dd y\,\dd x}{|x-y|^{n+s}},
 \quad 0<s<1.
\]

Several well-known indices recover classical integral-geometric functionals. Recall that
\begin{equation}\label{special-values}
 I_0(K)=\frac{\omega_{n-1}}{n\omega_n}S(K),
 \qquad
 I_1(K)=|K|,
 \qquad
 I_{n+1}(K)=\frac{n+1}{\omega_n}|K|^2.
\end{equation}
These identities are, respectively, Cauchy's projection formula, Crofton's formula for volume, and the Poincar\'e-Hadwiger formula, see \cites{Ren,Santalo}.
At the higher index treated in this paper, we show that
\begin{equation}\label{intro-inertia}
 I_{n+3}(K)
 =\frac{2(n+3)(n+2)}{n\omega_n}\,|K|\mathcal{I}(K),
\end{equation}
where $\mathcal{I}(K)$ is the polar moment of inertia of $K$ about its centroid $c_K$, i.e.,
\[
 c_K=\frac{1}{|K|}\int_Kx\,\dd x,
 \qquad
 \mathcal{I}(K)=\int_K|x-c_K|^2\,\dd x.
\]
Moreover, Zhang's dual kinematic formula \cite{Zhang99} gives, for $1\le m<n$,
\begin{equation}\label{intro-Zhang}
 \int_{\mathcal A(n,m)} |K\cap F|^2\,\dd F
 = \frac{m\omega_m}{n\omega_n} \int_K\int_K |x-y|^{m-n}\,\dd x\,\dd y
 = \frac{\omega_m}{m+1}I_{m+1}(K),
\end{equation}
where $\dd F$ is the invariant measure on the affine Grassmannian $\mathcal A(n,m)$, normalized compatibly with \eqref{def-Iq}.
% Together, these identities connect moments of random line sections with averages of affine-section volumes and show that the family $\{I_q\}$ interpolates among classical projection and volume functionals, affine-section moments, pair-distance energies, and inertia functionals.
Together, these identities show that the family $\{I_q\}$ links random line-section moments and affine-section averages, while interpolating among classical projection and volume functionals, pair-distance energies, and inertia functionals.

The first variation of $I_q$ gives rise to a family of geometric measures. More precisely, the $q$-th \textit{chord measure} $F_q(K,\cdot)$ is characterized by
\begin{equation}\label{intro-variation}
 \left.\frac{\dd}{\dd s}\right|_{s=0+} I_q(K+sL)
 =\int_{\mathbb S^{n-1}}h_L(u)\,\dd F_q(K,u),
\end{equation}
where $h_L$ is the support function of $L$. Lutwak, Xi, Yang, and Zhang \cite{LXYZ} introduced chord measures and solved the associated chord Minkowski problem: for $q>0$, a nonzero finite Borel measure $\mu$ on $\mathbb S^{n-1}$ is of the form
\[
 F_q(K,\cdot)=\mu
\]
for some $K\in\mathcal{K}^n$ if and only if $\mu$ is not concentrated on a closed hemisphere and has its centroid at the origin. 
They also introduced the associated $L_p$ chord measures and connected the resulting Minkowski problems with fully nonlinear Monge-Amp\`ere equations involving the nonlocal potentials $\widetilde V_{q-1}(K,z)$.

The existence theorem naturally leads to the corresponding uniqueness problem. In particular, one may ask whether
\[
 F_q(K,\cdot)=F_q(L,\cdot)
\]
implies that $K$ and $L$ differ only by a translation. This question was raised in \cites{LXYZ,HYZ}.
%and reiterated by Huang, Yang, and Zhang \cite{HYZ}.
As a natural counterpart to the uniqueness problem for chord measures, Huang, Yang, and Zhang \cite{HYZ} proposed studying Brunn-Minkowski inequalities for chord integrals.
% The survery \cite[Section~7.6]{HYZ} proposed studying Brunn-Minkowski inequalities for chord integrals as a counterpart to
% the uniqueness problem for chord measures.
For $K,L\in\mathcal{K}^n$, the natural Brunn-Minkowski inequality for $I_q$ is
\begin{equation}\label{BM-Iq-0}
 I_q((1-t)K+tL)^{\frac{1}{d}}
 \geq (1-t)I_q(K)^{\frac{1}{d}} +tI_q(L)^{\frac{1}{d}},\qquad 0\leq t\leq 1,
\end{equation}
where $d=n+q-1$ is the degree of homogeneity of $I_q$.
At $q=1$ and $q=n+1$, it follows from \eqref{special-values} and the classical Brunn-Minkowski inequality for volume. The case $q=0$ is related to the Brunn-Minkowski inequality for surface area, which is a consequence of the Aleksandrov-Fenchel inequality. 
In this paper, we determine the further range of $q$ for which \eqref{BM-Iq-0} holds and establish its equality cases.

Our first main result is the following.

\begin{thm}\label{thm-main}
Let $K,L\in\mathcal K^n$ and $0<t<1$. For every $q\in[1,n+2]\cup\{n+3\}$,
\begin{equation}\label{BM-Iq}
     I_q((1-t)K+tL)^{\frac{1}{n+q-1}}
 \geq (1-t)I_q(K)^{\frac{1}{n+q-1}} +tI_q(L)^{\frac{1}{n+q-1}}.
\end{equation}
Equality holds if and only if $K$ and $L$ are homothetic.
In dimension two, the same assertions hold for every $1\leq q\leq5$.
\end{thm}

\begin{rem}
After we have completed our paper, we note that in \cite{LLS26}, Liu-Lu-Sun proved the logarithmic Brunn-Minkowski inequality for chord integrals of origin-symmetric convex bodies for $1<q<n+1$, with equality precisely for dilates.
Their result implies \eqref{BM-Iq} in this symmetric setting, while Theorem~\ref{thm-main} applies to arbitrary convex bodies and a larger range of $q$.
\end{rem}

The proof uses the different representations of chord integrals developed above. 
By the covariogram representation \eqref{intro-cov-formula}, joint covariogram concavity and the weighted-marginal theorem of Cordero-Erausquin and Eskenazis \cite{CEE} give the range $1<q\leq n+1$.
The $X$-ray representation \eqref{intro-Xray},  together with convex lifts and linear Brunn-Minkowski inequalities, provides an independent proof for $1<q\leq2$. A symmetric lift of the covariogram extends the range to $n+1<q\leq n+2$ and gives the equality cases.
Finally, at the endpoint $q=n+3$, the inertia representation \eqref{intro-inertia} and Hadwiger's inertia inequality give the result.

Furthermore, we establish a transfer principle from symmetric dual quermassintegrals to chord integrals via convolution bodies. Notably, the equality characterization for the transferred chord inequality does not require an equality characterization for the assumed symmetric dual quermassintegral inequality.

\begin{thm}%[Dual-to-chord transfer]
\label{thm-intro-transfer}
Let $q>1$, set $\beta=q-1$, and let $0<t<1$. Assume that
\begin{equation}\label{BM-tildeV}
    \widetilde V_{\beta}((1-t)C+tD)^{\frac{1}{\beta}}
 \geq (1-t)\widetilde V_{\beta}(C)^{\frac{1}{\beta}} +t\widetilde V_{\beta}(D)^{\frac{1}{\beta}}
\end{equation}
holds for all origin-symmetric $C,D\in\mathcal K^n$. Then the chord Brunn-Minkowski inequality \eqref{BM-Iq} holds for all $K,L\in\mathcal K^n$.
Equality holds if and only if $K$ and $L$ are homothetic.
\end{thm}

As an application, combining the transfer principle with the results of Sadovsky and Zhang \cite{SZ25}, our previous work \cite{LWZ}, and Corollary~\ref{cor-dual-BM}, which establishes \eqref{BM-tildeV} for $n<\beta\leq n+1$, we obtain an alternative proof of the general-dimensional assertion of Theorem~\ref{thm-main}.
In the plane, incorporating Tao's result \cite{Tao26} also yields the two-dimensional assertion of Theorem~\ref{thm-main}.

For the case $0<q<1$, we obtain the following partial results.

\begin{thm}%[Fractional comparison and simplex inequality]
\label{thm:fractional-main}
Let $0<q<1$. Define
\begin{equation}\label{eq:frac-T-definition}
 T_q(K)
 =|K|^q\int_{S^{n-1}}b_K(u)^{1-q}\,d\sigma(u),
\end{equation}
where $ b_K(u)=|P_{u^\perp}K|_{n-1}$ denotes the brightness of $K$ in the direction $u$.

\begin{enumerate}
\item
For every $K\in\mathcal K^n$,
\begin{equation}\label{eq:frac-sharp-comparison}
 c_{n,q}T_q(K)\leq I_q(K)\leq T_q(K),
\end{equation}
where $c_{n,q} = \frac{n^q\Gamma(q+1)\Gamma(n)}{\Gamma(n+q)}$.
Equality on the left holds exactly for simplices.
%Equality in the lower bound holds if and only if $K$ is a simplex.

\item
For every $K,L\in\mathcal K^n$ and $0<t<1$,
\begin{equation}\label{eq:Tq-BM}
 \begin{split}
 T_q\big((1-t)K+tL\big)^{\frac{1}{n+q-1}}
 \geq
 (1-t)T_q(K)^{\frac{1}{n+q-1}}+t T_q(L)^{\frac{1}{n+q-1}}.
 \end{split}
\end{equation}
Consequently,
\begin{equation}\label{eq:frac-factor-BM}
 \begin{split}
 I_q\big((1-t)K+tL\big)^{\frac{1}{n+q-1}}
 \geq
 c_{n,q}^{\frac{1}{n+q-1}}
 \big((1-t)I_q(K)^{\frac{1}{n+q-1}}+tI_q(L)^{\frac{1}{n+q-1}}\big).
 \end{split}
\end{equation}
Moreover, if $S_0,S_1\in\mathcal K^n$ are simplices, then
\begin{equation}\label{eq:frac-simplex-BM}
 \begin{split}
 I_q\big((1-t)S_0+tS_1\big)^{\frac{1}{n+q-1}}
 \geq {}&
 (1-t)I_q(S_0)^{\frac{1}{n+q-1}}+tI_q(S_1)^{\frac{1}{n+q-1}}.
 \end{split}
\end{equation}
Equality holds if and only if $S_0$ and $S_1$ are homothetic.
\end{enumerate}
\end{thm}

\begin{rem}
 $T_q$ can be expressed in terms of a dual quermassintegral of the polar projection body $\Pi^*K=(\Pi K)^\circ$.
Indeed, $h_{\Pi K}=b_K$ and $\rho_{\Pi^*K}=b_K^{-1}$, so
\[
 T_q(K)
 =\frac{1}{\omega_n}|K|^q\, \widetilde V_{q-1}(\Pi^*K).
\]
\end{rem}

At the opposite end of the parameter range, the chord Brunn-Minkowski inequality fails beyond an explicit threshold.

\begin{thm}\label{thm-negative-main}
Let $\beta_*\in(2,3)$ be the constant defined by \eqref{threshold-definition}. Then \eqref{BM-Iq} fails for every $q>n+1+\beta_*.$
Numerically, $\beta_*\simeq 2.881467.$
\end{thm}

The construction is based on the rectangular perturbation method developed in our previous work \cite{LWZ} for dual quermassintegrals.
In the present setting, the relevant second variation is governed by the self-covariogram of a rectangle, whose triangular density leads to a different threshold. A one-dimensional reduction and a covariance estimate identify the corresponding sign change at $\beta_*$. In particular, \eqref{BM-Iq} fails for every $q\geq n+4$. The counterexamples can be chosen smooth, strictly convex, origin-symmetric, and unconditional.

As an application, the chord Brunn-Minkowski inequality and its $L_p$ counterpart yield uniqueness results for the chord Minkowski problems. Recall that
\[
 \dd F_{p,q}(K,\cdot)
 =h_K^{1-p}\,\dd F_q(K,\cdot).
\]
For $p>1$, the $L_p$ chord Brunn-Minkowski inequality holds for every $q\in[1,n+2]\cup\{n+3\}$, and for $1\leq q\leq 5$ in the plane, with equality precisely for dilates; see Theorem~\ref{thm-Firey}.

\begin{thm}\label{thm-chord-uniqueness}
Let $q\in[1,n+2]\cup\{n+3\}$. In the planar case, any $1\leq q\leq 5$ is admissible.
\begin{enumerate}
\item If $K,L\in\mathcal K^n$ and $F_q(K,\cdot)=F_q(L,\cdot)$, then $K,L$ differ by a translation.
\item If $K,L\in\mathcal K_o^n$, $p>1$, and $F_{p,q}(K,\cdot)=F_{p,q}(L,\cdot)$, then $K=L$ when $p\ne n+q-1$, and $K$ and $L$ are dilates when $p=n+q-1$.
\end{enumerate}
\end{thm}

% thm unique non-unique

Finally, we consider two extensions. The first concerns moments of affine $m$-sections,
\[
 I_{m,q}(K)
 =\int_{\mathcal A(n,m)}|K\cap F|_m^q\,\dd F.
\]
The whole-body lift yields a Brunn-Minkowski inequality for $1\leq q\leq1+1/m$, while Zhang's formula \eqref{intro-Zhang} gives $q=2$.
The second extension concerns mixed Riesz energies of the form
\[
 \mathcal E_\beta(K,L)
 =\int_K\int_L |x-y|^{\beta}\,\dd x\,\dd y, \qquad \beta>-n.
\]
The mixed covariogram $g_{K,L}(z)=|K\cap(L+z)|$ satisfies the same joint $1/n$-concavity property of $g_K=g_{K,K}$.
Along paths with even mixed covariograms, weighted marginals, symmetric lifts, and the level-set transfer yield joint Brunn-Minkowski inequalities for $\beta\in(-n,1]\cup\{2\}$, and for every $-2<\beta\leq2$ in the plane.

The paper is organized as follows. Section~\ref{sec:preliminaries} collects notation, identities, and normalization principles. Sections~\ref{sec:weighted-marginals}--\ref{sec:inertia-counterexamples} prove the general-dimensional part of Theorem~\ref{thm-main} using weighted marginals, linear Brunn-Minkowski lifts, and the inertia argument. Section~\ref{sec:inertia-counterexamples} also gives the rectangular counterexamples in Theorem~\ref{thm-negative-main}.
Section~\ref{sec:dual-to-chord} proves the transfer principle in Theorem~\ref{thm-intro-transfer} and the planar part of Theorem~\ref{thm-main}. 
Section~\ref{sec:fractional-simplices} proves the fractional comparison and simplex inequalities in Theorem~\ref{thm:fractional-main}. Section~\ref{sec:consequences} derives Minkowski inequalities and proves the uniqueness results in Theorem~\ref{thm-chord-uniqueness}. 
Section~\ref{sec:extensions} considers affine sections and mixed energies. 
Finally, Section~\ref{sec:problems} records the remaining problems.

\begin{ack}
    H. Li and Y. Zhang were supported by NSFC grant No. 12471047.
    Y. Wan was supported by Hong Kong RGC grants No. 14300623, No. 14304225, and No. 14308325.
\end{ack}

\section{Preliminaries}\label{sec:preliminaries}

\subsection{Convex bodies}\,

Let $\mathcal K^n$ denote the class of compact convex subsets of $\mathbb R^n$ with nonempty interior. Let $\mathcal K_o^n$ consist of those containing the origin in their interiors, and let $ \mathcal K_e^n$ denote the class of origin-symmetric convex bodies. The support function of $K\in\mathcal K^n$ is
\[
 h_K(u)=\max_{x\in K}\langle x,u\rangle,
 \qquad u\in\mathbb S^{n-1}.
\]
For $z\in K$, its radial function with respect to $z$ is
\[
 \rho_{K,z}(u)=\max\{r\geq0:z+ru\in K\}.
\]
We write $\rho_K=\rho_{K,0}$ when $0\in K$.
The covariogram function of $K$ is
\[
 g_K(z)
 =|K\cap(K+z)|
 =(\chi_K*\chi_{-K})(z),
 \qquad z\in\mathbb R^n,
\]
where $\chi_K$ denotes the indicator function of $K$. It is continuous and even, with
\[
 g_K(0)=|K|,
 \qquad
 \operatorname{supp}g_K=K-K.
\]
Moreover, $g_K>0$ on $\operatorname{int}(K-K)$, and $g_K^{1/n}$ is concave on $K-K$ by the Brunn-Minkowski inequality.

We denote Lebesgue measure by $|\cdot|$, with a dimensional subscript when needed, and write $\mathcal H^{n-1}$ for surface measure on $\partial K$. Let $\omega_j=|B^j|$. On $\mathbb S^{n-1}$, $\dd u$ denotes surface measure and $\dd\sigma$ denotes probability measure, so that
\[
 \dd u=n\omega_n\,\dd\sigma(u).
\]
The Grassmannian $G(n,m)$ carries probability Haar measure.
The invariant measure on the affine Grassmannian $\mathcal A(n,m)$ is normalized by
\[
 \int_{\mathcal A(n,m)}f(F)\,\dd F
 = \int_{G(n,m)}\int_{E^\perp} f(y+E)\,\dd y\,\dd\nu_m(E),
\]
where $\nu_m$ is probability Haar measure on $G(n,m)$.

\subsection{Chord integrals, potentials, and chord measures}\,

We use the normalization of chord integrals and chord measures in \cite{LXYZ}; see also \cite{GZ98,Zhang99,HYZ} for their geometric interpretations and integral representations.
For $q>0$, $I_q$ is positive, continuous in the Hausdorff metric, and strictly increasing under proper containment in $\mathcal K^n$.
%The last property follows from \eqref{intro-Xray}, since proper containment gives a strict chord-length inequality on an open set of lines.

For $z\in K$ and $r\in\mathbb R$, write
\[
 \widetilde V_r(K,z)
 =\frac1n \int_{\{u\in\mathbb S^{n-1}:\rho_{K,z}(u)>0\}}\rho_{K,z}(u)^r\,\dd u,
\]
whenever the integral is finite. If $z\in\operatorname{int}K$, the domain of integration is the whole sphere. This convention is particularly relevant at boundary points when $r=0$.

For $\beta>-n$, define the power-kernel potential
\begin{equation}\label{def-potential}
 U_{\beta,K}(x)
 = \int_K|x-y|^\beta\,\dd y, \qquad x\in\mathbb R^n.
\end{equation}
The kernel is locally integrable, so $U_{\beta,K}$ is finite and continuous. For $-n<\beta<0$, it is the Riesz potential of $\chi_K$ of order $\alpha=n+\beta$, without the usual normalization constant.
Polar coordinates give
\begin{equation}\label{potential-dual}
 \widetilde V_r(K,x)
 = \frac{r}{n}U_{r-n,K}(x), \qquad r>0,\quad x\in K.
\end{equation}
Thus the potential representation also extends $\widetilde V_r(K,\cdot)$ continuously to $\mathbb R^n$ for $r>0$; see \cite{LXYZ,CLL26}.
Integrating the potential over $K$ gives the self-energy
\begin{equation}\label{def-energy}
 \begin{split}
 \mathcal E_\beta(K)
 &= \int_K U_{\beta,K}(x)\,\dd x
 =\int_K\int_K|x-y|^\beta\,\dd x\,\dd y\\
 &= \int_{\mathbb R^n}g_K(z)|z|^\beta\,\dd z,
 \qquad \beta>-n.
 \end{split}
\end{equation}
The last identity follows from convolution. In particular, we have
\begin{equation}\label{energy-hom}
 \mathcal E_\beta(rK)=r^{2n+\beta}\mathcal E_\beta(K),
 \qquad
 \mathcal E_0(K)=|K|^2.
\end{equation}
For $n=3$, $\mathcal E_{-1}(K)$ is the Newtonian interaction integral of a body with uniform unit density.
Combining \eqref{potential-dual} with \eqref{intro-dual} yields
\begin{equation}\label{chord-energy}
 I_q(K)
 = \frac{q(q-1)}{n\omega_n}\mathcal E_{q-n-1}(K),
 \qquad q>1.
\end{equation}
% In Riesz-order notation, $\alpha=q-1$ and the homogeneity degree is $n+\alpha=n+q-1$. No symmetry assumption on $K$ is required.
For $0<q<1$, the corresponding interior potential is divergent; the finite exterior representation is \eqref{intro-fractional}.

Let $\nu_K$ denote the outer unit normal, defined $\mathcal H^{n-1}$-almost everywhere on $\partial K$. For $q>0$, the chord measure has the boundary representation
\begin{equation}\label{boundary-chord-measure}
 F_q(K,\eta)
 =\frac{2q}{\omega_n} \int_{\nu_K^{-1}(\eta)} \widetilde V_{q-1}(K,x)\,\dd\mathcal H^{n-1}(x),
\end{equation}
for every Borel set $\eta\subset\mathbb S^{n-1}$; see \cite{LXYZ}. The boundary integral defines a finite measure even when $0<q<1$. For $q>1$, \eqref{potential-dual} gives the equivalent potential formula
\[
 F_q(K,\eta)
 =\frac{2q(q-1)}{n\omega_n}\int_{\nu_K^{-1}(\eta)} U_{q-n-1,K}(x)\,\dd\mathcal H^{n-1}(x).
\]
In particular,
\[
 F_1(K,\cdot)=S(K,\cdot),
 \qquad
 F_{n+1}(K,\cdot)
 = \frac{2(n+1)}{\omega_n}|K|S(K,\cdot),
\]
where $S(K,\cdot)$ is the surface area measure.

The first-variation formula \eqref{intro-variation} implies
\begin{equation}\label{chord-euler}
 \int_{\mathbb S^{n-1}}h_K(u)\,\dd F_q(K,u)
 = (n+q-1)I_q(K).
\end{equation}
Chord measures are translation invariant, homogeneous of degree $n+q-2$, and satisfy
\[
 \int_{\mathbb S^{n-1}}u\,\dd F_q(K,u)=0,
\]
which follows by differentiating the translation invariance of $I_q$.

For $K\in\mathcal K_o^n$ and $p\in\mathbb R$, the $L_p$ chord measure is
\[
 \dd F_{p,q}(K,u) = h_K(u)^{1-p}\,\dd F_q(K,u).
\]
% It satisfies
% \[
%  F_{p,q}(rK,\cdot)
%  =
%  r^{n+q-p-1}F_{p,q}(K,\cdot).
% \]
Prescribing $F_{p,q}$ is the $L_p$ chord Minkowski problem.
The case $p=1$ is the ordinary chord Minkowski problem, whose existence theory was established in \cite{LXYZ}.
Further existence results were obtained in \cite{XYZ+Z} for $p>1$ and for even data when $0<p<1$, and in \cite{GXZ} for $0\leq p<1$ without symmetry assumptions.
Chen, Li, and Li \cite{CLL26} treated the supercritical range $p<-n-q+1$ for $3<q<n+1$ under suitable assumptions on the density. 
Nonuniqueness can occur for negative $p$; see \cite{Li24}.
In this paper, we study the uniqueness problem for $p\geq 1$.

\subsection{Homogeneity and normalization}\,

We collect several elementary principles used below. Let $\mathcal C$ be a class closed under Minkowski addition and positive dilations, such as $\mathcal K^n$ or $\mathcal K_o^n$.

\begin{prop}\label{prop-calculus}
Let $F:\mathcal C\to(0,\infty)$ be homogeneous of degree $d>0$.
\begin{enumerate}[(i)]
    \item If $F^a$ is concave for some $a>0$, then $a\leq1/d$ and $F^b$ is concave for every $0<b\leq a$. Moreover, $F^{1/d}$ is concave if and only if
\[
 F((1-t)K+tL)\geq1 \quad\text{whenever}\quad F(K)=F(L)=1.
\]
In particular, concavity of any positive power of $F$ implies concavity of  $F^{1/d}$.

    \item Suppose $\mathcal C\subset\mathcal K^n$, $d>1$, and the following linear inequality holds:
\begin{equation}\label{F-linear-BM}
     |K|=|L|
 \quad\Longrightarrow\quad
 F((1-t)K+tL)\geq(1-t)F(K)+tF(L).
\end{equation}
    Set $a=|K|^{1/n}$, $b=|L|^{1/n}$, and $c=(1-t)a+tb$. Then, for $0<t<1$,
\begin{equation}\label{strengthened-BM}
 F((1-t)K+tL)
 \geq c^{d-1} \left( (1-t)\frac{F(K)}{a^{d-1}} +t\frac{F(L)}{b^{d-1}} \right).
\end{equation}
    Consequently, $F^{1/d}$ is concave. Equality in the inequality for $F^{1/d}$ forces equality in the linear inequality \eqref{F-linear-BM}.

    \item If $G:\mathcal C\to(0,\infty)$ is homogeneous of degree $e>0$ and $F^{1/d},G^{1/e}$ are concave, then
\[
 (F^rG^s)^{\frac{1}{rd+se}}
\]
is concave for every $r,s>0$.

    \item Suppose $\mathcal C=\mathcal K_o^n$, $F$ is increasing under inclusion, and $F^{1/d}$ is concave. Then, for $p\geq1$ and $0< t<1$,
\begin{equation}\label{general-Firey}
 F((1-t)\cdot_pK+_pt\cdot_pL)^{\frac{p}{d}}
 \geq (1-t)F(K)^{\frac{p}{d}}+tF(L)^{\frac{p}{d}}.
\end{equation}
If $p>1$ and $F$ is strictly increasing, equality holds if and only if $K$ and $L$ are dilates.
\end{enumerate}
\end{prop}

\begin{proof}
For (i), concavity of $F^b$ follows by composing $F^a$ with the increasing concave function $s\mapsto s^{b/a}$. Restriction to a dilation ray $\{rK\}$ gives concavity of $r\mapsto r^{da}$ and hence $da\leq1$. 
To prove the normalization assertion, set $\tilde{a}=F(K)^{1/d}$, $\tilde{b}=F(L)^{1/d}$ and $\tilde{c}=(1-t)\tilde{a}+t\tilde{b}$. Under the change of variables $s=t\tilde{b}/\tilde{c}$, we have
\[
 \frac{(1-t)K+tL}{\tilde{c}}
 = (1-s)\frac{K}{\tilde{a}}+s\frac{L}{\tilde{b}}.
\]
Since $F(K/\tilde{a})=F(L/\tilde{b})=1$, the assumed superlevel-set condition gives $F((1-t)K+tL)\geq \tilde{c}^d$. Concavity of any positive power implies this condition.

For (ii), apply the equal-volume linear inequality to $K/a$ and $L/b$ with parameter $s=tb/c$. Homogeneity gives \eqref{strengthened-BM}. H\"older's inequality yields
\[
 \bigl((1-t)F(K)^{\frac{1}{d}}+tF(L)^{\frac{1}{d}}\bigr)^d
 \leq  c^{d-1} \left( (1-t)\frac{F(K)}{a^{d-1}} +t\frac{F(L)}{b^{d-1}} \right).
\]
This proves the concavity of $F^{1/d}$ and the equality implication.

For (iii), write
\[
 (F^rG^s)^{\frac{1}{rd+se} }
 = \bigl(F^{\frac{1}{d}}\bigr)^{\frac{rd}{rd+se} }\bigl(G^{\frac{1}{e}}\bigr)^{\frac{se}{rd+se}}
\]
and use the monotonicity and concavity of the weighted geometric mean.

Part (iv) is the $L_p$ transference principle of Zou and Xiong \cite{ZX18}; we include the argument for completeness.
Suppose $p>1$ and put $p'=p/(p-1)$. For $a,b\geq0$ with $a^{p'}+b^{p'}=1$, H\"older's inequality for support functions gives
\[
 a(1-t)^{1/p}K+bt^{1/p}L
 \subset (1-t)\cdot_pK+_pt\cdot_pL.
\]
Applying the monotonicity and concavity of $F^{1/d}$, we obtain
\[
 F((1-t)\cdot_pK+_pt\cdot_pL)^{\frac{1}{d}}
 \geq a(1-t)^{1/p}F(K)^{\frac{1}{d}} +bt^{1/p}F(L)^{\frac{1}{d}}.
\]
Maximizing over $a,b$ proves \eqref{general-Firey}.
If equality holds and $F$ is strictly increasing, the containment must be an equality. Equality in support-function H\"older then forces $h_L/h_K$ to be constant. Conversely, equality holds for dilates by homogeneity.
\end{proof}

\section{Weighted marginals and the range \texorpdfstring{$1< q\leq n+1$}{1 < q <= n+1}}\label{sec:weighted-marginals}

The covariogram representation \eqref{intro-cov-formula} expresses chord integrals as weighted integrals of $g_K$.
For $1<q\leq n+1$, the weight is a nonincreasing radial power, suggesting an application of the weighted-marginal theorem of Cordero-Erausquin and Eskenazis \cite{CEE}. We use the following form.

\begin{thm}[\cite{CEE}]%[Weighted marginals]
\label{thm-WM}
Let $\Omega\subset\mathbb R\times\mathbb R^n$ be convex, with origin-symmetric slices
\[
 \Omega_t=\{z:(t,z)\in\Omega\}.
\]
Suppose that $\Phi:\Omega\to[0,\infty)$ is concave and $\Phi(t,\cdot)$ is even for every $t$. Let $w:[0,\infty)\to\mathbb R$ be nondecreasing, with $u\mapsto w(e^u)$ convex. For $b>0$, the function
\[
 t\longmapsto
 \left(
 \int_{\Omega_t}\Phi(t,z)^b e^{-w(|z|)}\,\dd z
 \right)^{\frac{1}{n+b}}
\]
is concave whenever the integrals are finite.
\end{thm}

To apply this theorem along $K_t=(1-t)K+tL$, we need concavity jointly in $t$ and $z$, rather than concavity of each covariogram separately.

\begin{lem}\label{lem-cov-concavity}
Let $K,L\in\mathcal K^n$. The set
\[
 \Omega=\{(t,z):0\leq t\leq1,\ z\in K_t-K_t\}
\]
is convex, and the function
\[
 \Phi(t,z)=g_{K_t}(z)^{\frac{1}{n}}
\]
is concave on $\Omega$. Moreover, $\Phi(t,\cdot)$ is even for every $t$.
\end{lem}

\begin{proof}
For $(t,z)\in\Omega$, consider the centered intersection
\[
 C(t,z) = \{x\in\mathbb R^n:x\pm z/2\in K_t\}.
\]
This set is a translation of $K_t\cap(K_t+z)$, so $|C(t,z)|=g_{K_t}(z)$. 
Moreover, $ C(t,z)$ is nonempty if and only if $(t,z)\in \Omega$.

For $(t_i,z_i)\in\Omega$, $i=0,1$, and $0\leq s\leq1$, the identity $ (1-s)K_{t_0}+sK_{t_1} = K_{(1-s)t_0+st_1}$ gives
\[
 (1-s)C(t_0,z_0)+sC(t_1,z_1)
 \subset
 C\bigl((1-s)t_0+st_1,(1-s)z_0+sz_1\bigr).
\]
Hence $\Omega$ is convex. The classical Brunn-Minkowski inequality then implies that $\Phi$ is concave. Finally, the symmetry $C(t,z)=C(t,-z)$ gives evenness.
\end{proof}

To apply Theorem~\ref{thm-WM} to the covariogram representation, Lemma~\ref{lem-cov-concavity} is sufficient to prove the chord Brunn-Minkowski inequality for $1< q\leq n+1$. However, the weighted-marginal argument alone does not provide an equality characterization. To recover rigidity, we replace the interpolated covariogram by the jointly concave envelope generated by the two endpoint covariograms.

For $K_t=(1-t)K+tL$, set
\[
 D_K=K-K,\qquad D_L=L-L,\qquad
 D_t=(1-t)D_K+tD_L.
\]
For $z\in D_t$, define
\begin{equation}\label{eq:psi-definition}
 \psi_t(z)
 =\sup_{\substack{x\in D_K,\ y\in D_L\\
                  z=(1-t)x+ty}}
 \left[ (1-t)g_K(x)^{\frac{1}{n}}+t g_L(y)^{\frac{1}{n}} \right].
\end{equation}
By compactness, the supremum is attained. Lemma~\ref{lem-cov-concavity} then immediately yields
\begin{equation}\label{eq:psi-covariogram-bound}
 \psi_t(z)\leq g_{K_t}(z)^{\frac{1}{n}}.
\end{equation}
Furthermore, we have
\begin{lem}\label{lem-psi-concavity}
The function $(t,z)\mapsto\psi_t(z)$ is concave on $ \{(t,z):0\leq t\leq1,\ z\in D_t\}$, and $\psi_t$ is even for every $t$. Moreover,
\[
 \psi_t(0)= (1-t)|K|^{\frac{1}{n}}+t|L|^{\frac{1}{n}}.
\]
\end{lem}

\begin{proof}
Concavity follows directly from the definition. Indeed, let
\[
 z_i=(1-t_i)x_i+t_i y_i,\qquad 
 x_i\in D_K,\quad y_i\in D_L,\quad i=0,1,
\]
and set $t_s=(1-s)t_0+st_1$ for $0\leq s\leq 1$.
For $0<t_s<1$, define
\[
 x_s=
\frac{(1-s)(1-t_0)x_0+s(1-t_1)x_1}{1-t_s},
 \qquad
 y_s=\frac{(1-s)t_0y_0+st_1y_1}{t_s}.
\]
By the convexity of $D_K$ and $D_L$, we have $x_s\in D_K$ and $y_s\in D_L$, and
\[
 (1-s)z_0+sz_1=(1-t_s)x_s+t_sy_s.
\]
The concavity of $g_K^{1/n}$ and $g_L^{1/n}$ therefore gives
\[
\begin{aligned}
 &(1-t_s)g_K(x_s)^{\frac{1}{n}}+t_sg_L(y_s)^{\frac{1}{n}}\\
\geq\, &(1-s)\bigl((1-t_0)g_K(x_0)^{\frac{1}{n}} +t_0g_L(y_0)^{\frac{1}{n}}\bigr)
 +s\bigl((1-t_1)g_K(x_1)^{\frac{1}{n}} +t_1g_L(y_1)^{\frac{1}{n}}\bigr).
\end{aligned}
\]
Taking the supremum over admissible decompositions proves the joint concavity of $(t,z)\mapsto\psi_t(z)$.
Evenness follows from the evenness of $g_K$ and $g_L$.

Since $g_K$ and $g_L$ attain their maxima at the origin, the supremum in \eqref{eq:psi-definition} at $z=0$ is attained at $x=y=0$, proving the last identity.
\end{proof}

\begin{thm}\label{thm-negative-energy}
Let $K,L\in\mathcal K^n$ and $0< t<1$. For $-n<\beta\leq0$,
\begin{equation}\label{BM-energy-negative}
 \mathcal E_\beta(K_t)^{\frac{1}{2n+\beta}}
 \geq (1-t)\mathcal E_\beta(K)^{\frac{1}{2n+\beta}} +t\mathcal E_\beta(L)^{\frac{1}{2n+\beta}}.
\end{equation}
Consequently, the chord Brunn-Minkowski inequality \eqref{BM-Iq} holds for $1< q\leq n+1$. Equality holds if and only if $K$ and $L$ are homothetic.
\end{thm}

\begin{proof}
For $\beta=0$, the assertion and its equality characterization follow from $\mathcal E_0(K)=|K|^2$ and the classical Brunn-Minkowski inequality.

Assume that $-n<\beta<0$. For $\varepsilon>0$, set
\[
 W_\varepsilon(z)
 = (\varepsilon^2+|z|^2)^{\frac{\beta}{2}},
 \qquad
 w_\varepsilon(r)
 = -\frac{\beta}{2}\log(\varepsilon^2+r^2).
\]
Then $w_\varepsilon$ is nondecreasing and $u\mapsto w_\varepsilon(e^u)$ is convex. Applying Theorem~\ref{thm-WM} to $\psi_t$ shows that
\[
 t\longmapsto
 \left(
 \int_{D_t}\psi_t(z)^nW_\varepsilon(z)\,\dd z
 \right)^{\frac{1}{2n}}
\]
is concave.
Since $\psi_0=g_K^{1/n}$ and $\psi_1=g_L^{1/n}$, monotone convergence as $\varepsilon\to 0^+$ together with \eqref{eq:psi-covariogram-bound} yields
\begin{align*}
    (1-t)\mathcal E_\beta(K)^{\frac{1}{2n}} +t\mathcal E_\beta(L)^{\frac{1}{2n}}
    & \leq  \left(\int_{D_t}\psi_t(z)^n|z|^\beta\,\dd z\right)^{\frac{1}{2n}}\\
    & \leq \left(\int_{D_t}g_{K_t}(z)|z|^\beta\,\dd z\right)^{\frac{1}{2n}}
    =\mathcal E_\beta(K_t)^{\frac{1}{2n}}.
\end{align*}
Thus, Proposition~\ref{prop-calculus}(i) gives \eqref{BM-energy-negative} and then \eqref{BM-Iq}.

It remains to consider equality. By the normalization in Proposition~\ref{prop-calculus}(i), it suffices to consider $\mathcal E_\beta(K)=\mathcal E_\beta(L)=1$.
If $K$ and $L$ are not homothetic, the classical Brunn-Minkowski inequality gives
\[
 g_{K_t}(0)^{\frac{1}{n}} = |K_t|^{\frac{1}{n}}
 > (1-t)|K|^{\frac{1}{n}}+t|L|^{\frac{1}{n}} =\psi_t(0).
\]
By continuity, this strict inequality persists on a neighborhood of the origin. Hence
\[
 \mathcal E_\beta(K_t)
 >\int_{D_t}\psi_t(z)^n|z|^\beta\,\dd z\geq1,
\]
which excludes equality. The converse follows from translation invariance and homogeneity.
\end{proof}

\section{Linear Brunn-Minkowski inequalities and convex lifts}
\label{sec:linear-lifts}

Our second approach combines convex lifts with linear Brunn-Minkowski inequalities. Based on a method developed in our work on dual quermassintegrals, it applies to geometric functionals with suitable integral representations and yields Brunn-Minkowski inequalities on corresponding parameter ranges.

For chord integrals, we use two such constructions. The $X$-ray representation \eqref{intro-Xray} gives the range $1<q\leq2$, while a lift of the covariogram gives $n+1<q\leq n+2$.
In both cases, we first establish a linear inequality under an equal-volume constraint and then apply Proposition~\ref{prop-calculus}(ii). The linear form allows us to average the directional inequalities before taking the natural root.

We use the following projection form of the linear Brunn-Minkowski inequality, due to Bonnesen, together with a consequence of its equality characterization; see \cite{HCYN,SGYN,BS}.

\begin{thm}[Linear Brunn-Minkowski inequality]
\label{thm-Bonnesen}
Let $A,B\subset\mathbb R^{n+1}$ be convex bodies, and let $\pi$ be orthogonal projection onto a hyperplane. If $|\pi A|=|\pi B|$, then
\[
 |(1-t)A+tB| \geq  (1-t)|A|+t|B|,
 \qquad 0\leq t\leq1.
\]
If equality holds for $0<t<1$, the projected bodies $\pi A$ and $\pi B$ are translates.
\end{thm}

\subsection{The range \texorpdfstring{$1<q\leq2$}{1 < q <= 2}}\,

For a fixed direction $u\in\mathbb S^{n-1}$, Brunn's section theorem implies that the chord-length function $X_K(\cdot,u)$ is concave on $P_{u^\perp}K$.
Moreover, section inclusion gives
\begin{equation}\label{chord-profile}
 X_{K_t}((1-t)y+tz,u)
 \geq (1-t)X_K(y,u)+tX_L(z,u)
\end{equation}
for $y\in P_{u^\perp}K$ and $z\in P_{u^\perp}L$.

\begin{thm}\label{thm-low}
For $1<q\leq2$, the chord Brunn-Minkowski inequality \eqref{BM-Iq} holds. Equality holds if and only if $K$ and $L$ are homothetic. Moreover, \eqref{strengthened-BM} holds with $F=I_q$ and $d=n+q-1$.
\end{thm}

\begin{proof}
Set $\alpha=q-1\in(0,1]$. For each $u\in\mathbb S^{n-1}$, define the convex lift
\[
 \mathcal L_{u,q}(K)
 = \left\{ (x,r)\in\mathbb R^n\times\mathbb R: x\in K,\quad 0\leq r\leq X_K(P_{u^\perp}x,u)^\alpha \right\}.
\]
Its convexity follows from concavity of the chord-length function and of $s\mapsto s^\alpha$. It follows from \eqref{chord-profile} that
\[
 (1-t)\mathcal L_{u,q}(K)+t\mathcal L_{u,q}(L)
 \subset \mathcal L_{u,q}(K_t).
\]
Its projection onto $\mathbb R^n$ is $K$. Moreover, we have
\[
 |\mathcal L_{u,q}(K)|
 =\int_{u^\perp}X_K(y,u)^q\,\dd y.
\]

Suppose first that $|K|=|L|$. Theorem~\ref{thm-Bonnesen} therefore yields
\begin{equation}\label{ineq-directional}
     \int_{u^\perp}X_{K_t}(y,u)^q\,\dd y
 \geq (1-t)\int_{u^\perp}X_K(y,u)^q\,\dd y +t\int_{u^\perp}X_L(y,u)^q\,\dd y
\end{equation}
for any $u\in\mathbb{S}^{n-1}$. 
By \eqref{intro-Xray}, integrating over $u$ proves the equal-volume linear inequality for $I_q$. Since its degree is $n+q-1$,
Proposition~\ref{prop-calculus}(ii) gives both \eqref{BM-Iq} and \eqref{strengthened-BM} for $I_q$.

Equality in \eqref{BM-Iq} implies equality in the directional lift inequalities \eqref{ineq-directional} for almost every direction after volume normalization.
Theorem~\ref{thm-Bonnesen} then shows that the normalized bodies are translates, so $K$ and $L$ are homothetic. The converse follows from translation invariance and homogeneity.
\end{proof}

% \begin{rem}
% The linear directional inequalities must be integrated before taking
% the natural root. Averaging directional root inequalities would
% require an $L^p$ triangle inequality in the opposite direction.
% \end{rem}

\subsection{Projected positive moments}\,

The second construction uses the following consequence of the same lifting argument.

\begin{lem}\label{lem-projected}
Let $P:\mathbb R^N\to\mathbb R^k$ be a nonzero linear map and $0<\alpha\leq1$. For origin-symmetric convex bodies, define
\[
 M_{P,\alpha}(A)=\int_A|Px|^\alpha\,\dd x.
\]
If $|A|=|B|$, then
\[
 M_{P,\alpha}((1-t)A+tB)
 \geq (1-t)M_{P,\alpha}(A)+tM_{P,\alpha}(B),\qquad 0<t<1.
\]
Consequently, $M_{P,\alpha}^{1/(N+\alpha)}$ is concave.
Equality in the linear inequality holds if and only if $A=B$; equality in the natural-root inequality holds if and only if $A$ and $B$ are dilates.
% Equality holds if and only if $A$ and $B$ are dilates.
\end{lem}

\begin{proof}
For a nonzero linear functional $\ell$, consider
\[
 D_\ell(A)
 = \left\{ (x,r): x\in A,\quad \ell(x)\geq0,\quad 0\leq r\leq\ell(x)^\alpha \right\}.
\]
This is a convex body whose projection onto $\mathbb R^N$ is $ A_\ell^+=A\cap\{\ell\geq0\}$.
It follows from the origin symmetry of $A$ that
\[
 |A_\ell^+|=\frac{|A|}{2},
 \qquad
 |D_\ell(A)|
 = \frac12\int_A|\ell(x)|^\alpha\,\dd x.
\]
Since $s\mapsto s^{\alpha}$ is concave when $0<\alpha\leq 1$, we obtain
\[
 (1-t)D_\ell(A)+tD_\ell(B)
 \subset D_\ell((1-t)A+tB).
\]
Theorem~\ref{thm-Bonnesen} thus proves the equal-volume linear inequality for the directional moment.

Apply this inequality with $\ell(x)=\langle Px,e\rangle$ and integrate over $e\in\mathbb S^{k-1}$. The identity
\[
 \int_{\mathbb S^{k-1}}
 |\langle Px,e\rangle|^\alpha\,\dd\sigma(e)
 = c_{k,\alpha}|Px|^\alpha,
 \qquad c_{k,\alpha}>0,
\]
gives the asserted linear inequality. 
Since $M_{P,\alpha}$ has degree $N+\alpha$, Proposition~\ref{prop-calculus}(ii) gives the concavity of $M_{P,\alpha}^{1/(N+\alpha)}$.

After volume normalization, equality and Theorem~\ref{thm-Bonnesen} imply that $A_\ell^+$ and $B_\ell^+$ are translates for almost every nonzero cutting functional $\ell$.
The translation lies in $\ker\ell$ and vanishes by the central symmetry of their sections there. Hence $A_\ell^+=B_\ell^+$, and origin symmetry gives $A=B$. Without volume normalization, this argument yields that $A$ and $B$ are dilates. The converse follows from homogeneity.
\end{proof}

As a direct application, the same lifting argument yields the symmetric Brunn-Minkowski inequality for dual quermassintegrals in the range $n<q\leq n+1$. Indeed, for $K\in\mathcal K_o^n$,
\[
 \widetilde V_q(K)
 = \frac1n\int_{\mathbb S^{n-1}}\rho_K(u)^q\,\dd u
 = \frac qn\int_K|x|^{q-n}\,\dd x.
\]
Taking $N=k=n$, $P=\operatorname{Id}$, and
$\alpha=q-n\in(0,1]$ in Lemma~\ref{lem-projected} gives the following.

\begin{cor}\label{cor-dual-BM}
Let $K,L\in\mathcal K_o^n$ be origin-symmetric. If $n<q\leq n+1$, then
\begin{equation}\label{dual-BM}
 \widetilde V_q((1-t)K+tL)^{\frac{1}{q}}
 \geq (1-t)\widetilde V_q(K)^{\frac{1}{q}} +t\widetilde V_q(L)^{\frac{1}{q}},
 \qquad 0< t< 1.
\end{equation}
Equality holds if and only if $K$ and $L$ are dilates.
\end{cor}

\subsection{The range \texorpdfstring{$n+1<q\leq n+2$}{n+1 < q <= n+2}}\,

Above $n+1$, the height in the $X$-ray lift is no longer a concave power. We instead apply Lemma~\ref{lem-projected} to a symmetric convex lift of the covariogram.
Define
\begin{equation}\label{def-cov-lift}
 \mathcal B(K)
 =\left\{(z,w)\in\mathbb R^n\times\mathbb R^n: z\in K-K,\quad |w|\leq g_K(z)^{\frac{1}{n}} \right\}.
\end{equation}
By Lemma~\ref{lem-cov-concavity}, this is an origin-symmetric convex body and
\begin{equation}\label{lift-inclusion}
 (1-t)\mathcal B(K)+t\mathcal B(L)
 \subset \mathcal B(K_t).
\end{equation}
Moreover, we have
\begin{align}
 \mathcal B(rK)&=r\mathcal B(K),\label{lift-scaling}\\
 |\mathcal B(K)|&=\omega_n|K|^2,\label{lift-volume}\\
 \int_{\mathcal B(K)}|z|^\alpha\,\dd z\,\dd w
 &=\omega_n\mathcal E_\alpha(K).\label{lift-moment}
\end{align}
% Thus the symmetry required by the projected-moment lemma comes
% from the covariogram rather than from the original body.

\begin{thm}\label{thm-positive}
Let $K,L\in\mathcal{K}^n$ and $0<t<1$. For $0<\alpha\leq1$,
\begin{equation}\label{BM-energy-positive}
 \mathcal E_\alpha(K_t)^{\frac{1}{2n+\alpha}}
 \geq (1-t)\mathcal E_\alpha(K)^{\frac{1}{2n+\alpha}} +t\mathcal E_\alpha(L)^{\frac{1}{2n+\alpha}}.
\end{equation}
Consequently, \eqref{BM-Iq} holds for $n+1<q\leq n+2$.
In both inequalities, equality holds if and only if $K$ and $L$ are homothetic.
\end{thm}

\begin{proof}
Suppose first that $|K|=|L|$. By \eqref{lift-volume}, the origin-symmetric convex bodies $\mathcal B(K)$ and $\mathcal B(L)$ have equal volume. 
Applying Lemma~\ref{lem-projected} with $P(z,w)=z$, by \eqref{lift-inclusion}, we obtain
\[
 \begin{split}
 \omega_n\mathcal E_\alpha(K_t)
 &=M_{P,\alpha}(\mathcal B(K_t))\\
 &\geq M_{P,\alpha}((1-t)\mathcal B(K)+t\mathcal B(L))\\
 &\geq (1-t)M_{P,\alpha}(\mathcal B(K))  +tM_{P,\alpha}(\mathcal B(L))\\
 &=\omega_n\bigl((1-t)\mathcal E_\alpha(K) +t\mathcal E_\alpha(L)\bigr).
 \end{split}
\]
Since $\mathcal E_\alpha$ has degree $2n+\alpha$, Proposition~\ref{prop-calculus}(ii) proves \eqref{BM-energy-positive}.

For equality, normalize $K,L$ so that their energies $\mathcal E_\alpha$ are equal.
Equality in the lifted inequality for $M_{P,\alpha}^{1/(2n+\alpha)}$ and Lemma~\ref{lem-projected} then imply that $\mathcal B(K)$ and $\mathcal B(L)$ are dilates.
Their equal projected moments force them to coincide.
By \eqref{lift-inclusion}, we obtain
\[
 \mathcal B(K_t)=\mathcal B(K)=\mathcal B(L).
\]
By \eqref{lift-volume}, $|K_t|=|K|=|L|$, so equality in the classical Brunn-Minkowski inequality makes the normalized bodies translates. Without normalization, we conclude that $K$ and $L$ are homothetic.
The converse follows from translation invariance and homogeneity.

Finally, taking $\alpha=q-n-1\in(0,1]$ in \eqref{chord-energy} gives the chord inequality \eqref{BM-Iq} and its equality characterization.
\end{proof}

\section{The inertia endpoint and counterexamples}\label{sec:inertia-counterexamples}

\subsection{The endpoint \texorpdfstring{$q=n+3$}{q=n+3}}\,

Let
\[
 c_K=\frac1{|K|}\int_Kx\,\dd x,
 \qquad
 \mathcal I(K)=\int_K|x-c_K|^2\,\dd x
\]
denote the centroid and the polar moment of inertia of $K$ about its centroid. Hadwiger's inertia inequality \cite{Hadwiger} states that
\begin{equation}\label{BM-inertia}
 \mathcal I(K_t)^{\frac{1}{n+2}}
 \geq (1-t)\mathcal I(K)^{\frac{1}{n+2}} +t\mathcal I(L)^{\frac{1}{n+2}}.
\end{equation}
Equality holds precisely for homothetic bodies; see also \cite[Theorem~5.1]{LWZ} for the inequality and its equality characterization. The connection with chord integrals is given by the following identity.

\begin{lem}\label{lem-inertia-identity}
For every $K\in\mathcal K^n$,
\begin{equation}\label{Iq-inertia}
 I_{n+3}(K)
 = \frac{2(n+3)(n+2)}{n\omega_n} |K|\mathcal I(K).
\end{equation}
\end{lem}

\begin{proof}
Since $\int_K(x-c_K)\,\dd x=0$, we obtain
\begin{align*}
    \mathcal E_2(K)
 &=\int_K\int_K\bigl(|x-c_K|^2+|y-c_K|^2 -2\langle x-c_K,y-c_K\rangle\bigr) \,\dd x\,\dd y\\
 &= 2|K|\int_K|x-c_K|^2\,\dd x =2|K|\mathcal I(K).
\end{align*}
The identity follows from \eqref{chord-energy}.
\end{proof}

Combining with \eqref{BM-inertia}, the classical Brunn-Minkowski inequality, and Proposition~\ref{prop-calculus}(iii), we obtain the following endpoint result. The equality characterization follows directly from the classical Brunn-Minkowski inequality.

\begin{thm}\label{prop-inertia}
The chord Brunn-Minkowski inequality \eqref{BM-Iq} holds at $q=n+3$. Equality holds if and only if $K$ and $L$ are homothetic.
\end{thm}

Together with Theorems~\ref{thm-negative-energy}, \ref{thm-low}, and \ref{thm-positive}, this proves the  general-dimensional part of Theorem~\ref{thm-main}; its planar extension is proved in Section \ref{sec:dual-to-chord}.

\subsection{Proof of Theorem~\ref{thm-negative-main}}\, %Rectangular perturbations

To construct counterexamples to the chord Brunn-Minkowski inequality, we adapt the rectangular perturbations from our earlier work \cite{LWZ} on dual quermassintegrals. Here the second variation is determined by the triangular density of the rectangle's self-covariogram.

For $|t|<1$, set
\[
 R_t=[-(1+t),1+t]\times[-(1-t),1-t],
 \qquad
 E_\beta(t)=\mathcal E_\beta(R_t),
 \quad \beta>0.
\]
The rectangles $R_t$ and $R_{-t}$ are rotations of one another, and
\[
 \frac{R_t+R_{-t}}2=R_0.
\]
Thus $E_\beta(t)>E_\beta(0)$ contradicts midpoint concavity of every positive power of $\mathcal E_\beta$.

Since the overlap lengths of two translated intervals are triangular, the covariogram of $R_t$ is given by
\[
 g_{R_t}(z)
 =\bigl(2(1+t)-|z_1|\bigr)_+ \bigl(2(1-t)-|z_2|\bigr)_+,
\]
where $a_+=\max\{a,0\}$. Substituting $z_1=2(1+t)x$ and $z_2=2(1-t)y$ into the covariogram formula \eqref{def-energy}, and writing $w(x,y)=(1-|x|)(1-|y|)$, we have
\begin{equation}\label{rectangle-energy}
 \begin{split}
 E_\beta(t)
 =2^{\beta+4}(1-t^2)^2
 \int_{[-1,1]^2} w(x,y)  \bigl((1+t)^2x^2+(1-t)^2y^2\bigr)^{\frac{\beta}{2}} \,\dd x\,\dd y.
 \end{split}
\end{equation}
% Differentiation under the integral is justified locally in $t$
% by an integrable multiple of $(x^2+y^2)^{\beta/2}$. 
By symmetry, $E_\beta'(0)=0$, and direct calculation yields
\begin{equation}\label{rect-second}
 \frac{E_\beta''(0)}{E_\beta(0)}
 = H(\beta)
 := \beta-4+\beta(\beta-2)A(\beta),
\end{equation}
where
\[
 A(\beta)
 = \frac{\int_{[0,1]^2}
 (1-x)(1-y)(x^2+y^2)^{\frac{\beta}{2}}
 \left(\frac{x^2-y^2}{x^2+y^2}\right)^2
 \,\dd x\,\dd y
 }{\int_{[0,1]^2}
 (1-x)(1-y)(x^2+y^2)^{\frac{\beta}{2}}
 \,\dd x\,\dd y }.
\]
% In particular, the squared volume factor $(1-t^2)^2$ in
% \eqref{rectangle-energy} distinguishes this variation from
% that of a radial moment.
Splitting $[0,1]^2$ along its diagonal and substituting $y=ux$,
we use the identity
\[
 \int_0^1x^{\beta+1}(1-x)(1-ux)\,\dd x
 = \frac{\beta+4-(\beta+2)u}{(\beta+2)(\beta+3)(\beta+4)}
\]
to obtain the following representation
\[
 A(\beta)
 = \frac{\int_0^1W_\beta(u)a(u)\,\dd u}{\int_0^1W_\beta(u)\,\dd u},
\]
where $W_\beta(u)$ and $a(u)$ are defined by
\begin{align*}
     W_\beta(u)
 = [\beta+4-(\beta+2)u](1+u^2)^{\frac{\beta}{2}},
 \qquad
 a(u)=\left(\frac{1-u^2}{1+u^2}\right)^2.
\end{align*}

The next lemma shows that there is a unique $\beta_*\in(2,3)$ satisfying
\begin{equation}\label{threshold-definition}
 H(\beta_*)=0.
\end{equation}

\begin{lem}\label{lem-threshold}
The function $H$ is strictly increasing on $[2,4]$ and satisfies $H(2)<0<H(3)$. Moreover, $H(\beta)>0$ for $\beta\geq4$.
Consequently, $H(\beta)>0$ for every $\beta>\beta_*$.
\end{lem}

\begin{proof}
Set $Z_\beta=\int_0^1W_\beta(u)\,\dd u$. Differentiating $A(\beta)$ gives
\[
 A'(\beta)
 =\frac1{Z_\beta}\int_0^1W_\beta(u)\bigl(a(u)-A(\beta)\bigr) \bigl(\psi_\beta(u)-\overline\psi_\beta\bigr)\,\dd u,
\]
where
\[
 \psi_\beta(u)
 = \partial_\beta\log W_\beta(u)
 = \frac{1-u}{\beta+4-(\beta+2)u} +\frac12\log(1+u^2),
 \qquad
 \overline\psi_\beta
 =\frac1{Z_\beta}\int_0^1W_\beta\psi_\beta\,\dd u.
\]
For $2\leq\beta\leq4$, we have
\[
 0\leq a\leq1,
 \qquad
 0\leq\psi_\beta(u)
 \leq\frac12(1-u+u^2)
 \leq\frac12.
\]
Since $\int_0^1W_\beta(a-A(\beta))\,\dd u=\int_0^1
 W_\beta(\psi_\beta-\overline\psi_\beta)\,\dd u=0$, we obtain
\[
 \begin{aligned}
 \frac1{Z_\beta}\int_0^1W_\beta(a-A(\beta))^2\,\dd u
 &\leq
 \frac1{Z_\beta}\int_0^1W_\beta(a-\tfrac12)^2\,\dd u
 \leq\frac14,\\
 \frac1{Z_\beta}\int_0^1
 W_\beta(\psi_\beta-\overline\psi_\beta)^2\,\dd u
 &\leq
 \frac1{Z_\beta}\int_0^1W_\beta(\psi_\beta-\tfrac14)^2\,\dd u
 \leq\frac1{16}.
 \end{aligned}
\]
Cauchy-Schwarz inequality therefore yields $|A'(\beta)|\leq1/8$. Hence, for any $2\leq \beta\leq 4$,
\[
 H'(\beta)
 \geq 1-\frac{\beta(\beta-2)}8+(2\beta-2)A(\beta)>0.
\]

Clearly, $H(2)=-2$. To verify $H(3)>0$, a direct computation gives
\[
 \begin{split}
 &\int_{[0,1]^2}
 (1-x)(1-y)(x^2+y^2)^{\frac{3}{2}}\Bigl(-1+3\bigl(\frac{x^2-y^2}{x^2+y^2}\bigr)^2\Bigr)\,\dd x\,\dd y\\
 % &\qquad=
 % \frac2{105}\int_0^1
 % \frac{(7-5u)(1-4u^2+u^4)}{\sqrt{1+u^2}}\,\dd u\\
 =\,& \frac{9\operatorname{arsinh}(1)-11\sqrt2+8}{20}>0.
 \end{split}
\]
% Indeed, Jensen's inequality gives
% \[
%  \operatorname{arsinh}(1)
%  =
%  \int_0^1(1+u^2)^{-1/2}\,\dd u
%  \geq\frac{\sqrt3}{2},
% \]
% and $8+9\sqrt3/2>11\sqrt2$.
Finally, $A(\beta)>0$ implies $H(\beta)>0$ for $\beta\geq4$.
\end{proof}

\begin{proof}[Proof of Theorem~\ref{thm-negative-main}]
Fix $\beta>\beta_*$. By Lemma~\ref{lem-threshold}, $E_\beta''(0)>0$, so
\[
 E_\beta(t)>E_\beta(0)
\]
for sufficiently small $t\neq0$. This gives the planar counterexample.

For $n>2$, take the thin boxes
\[
 K_{\varepsilon,t}
 =R_t\times[-\varepsilon,\varepsilon]^{n-2}.
\]
Their midpoint is $K_{\varepsilon,0}$, and the two endpoints are rotations of one another. Rescaling the thin coordinates
and applying dominated convergence, we get
\begin{align*}
    \mathcal E_\beta(K_{\varepsilon,t})=(2\varepsilon)^{2(n-2)}E_\beta(t)+o(\varepsilon^{2(n-2)}),\quad \text{as }\varepsilon\to 0^+.
\end{align*}
The strict planar gap therefore persists for sufficiently small $\varepsilon$. Taking $q=n+1+\beta$ in \eqref{chord-energy} proves failure of \eqref{BM-Iq}.

Finally, standard smoothing preserving the symmetries and rotation relation retains the strict midpoint gap by continuity.
Thus the counterexamples may be chosen smooth, strictly convex, origin-symmetric, and unconditional.
The approximation $\beta_*\simeq2.881467$ is obtained by direct numerical evaluation of \eqref{threshold-definition}.
\end{proof}

\section{From symmetric dual quermassintegrals to chord integrals}
\label{sec:dual-to-chord}

The preceding results suggest a correspondence between the Brunn-Minkowski inequality for $I_q$ and that for the dual quermassintegral $\widetilde V_{q-1}$.
This correspondence is made precise by the covariogram representation \eqref{intro-cov-formula}.
The key observation is that the positive superlevel sets of a covariogram are origin-symmetric convex bodies, even when the original body is not origin-symmetric.
These sets are usually called \textit{convolution bodies}. The notion is due to Kiener, who observed that these superlevel sets are convex; see \cite[Note~9.3]{Gardner} and the references therein.

For $K\in\mathcal K^n$ and $s\in(0,|K|^{1/n})$, define the convolution body
\begin{equation}\label{eq:convolution-body}
 K(s)
 = \left\{ z\in K-K:g_K(z)^{1/n}\geq s \right\}
 = \left\{ z\in\mathbb R^n:g_K(z)\geq s^n \right\}.
\end{equation}
Since $g_K$ is even and $g_K^{1/n}$ is concave on $K-K$, $K(s)$ is an origin-symmetric convex body. We use only positive levels; at level $s=0$, the natural set is the difference body $D_K=K-K$.

For $\beta>0$, let $\dd\mu_\beta(z)=|z|^{\beta-n}\,\dd z.$
For $C\in\mathcal K_o^n$, polar coordinates give
\begin{equation}\label{eq:mu-dual-volume}
 \mu_\beta(C)
 =\frac1\beta \int_{\mathbb S^{n-1}}\rho_C(u)^\beta\,\dd u
 =\frac n\beta\,\widetilde V_\beta(C).
\end{equation}
Consequently, the Brunn-Minkowski inequality for $\widetilde V_\beta$ is equivalent to that for $\mu_\beta$.

We first record the level-set inclusion that follows from the joint covariogram concavity proved in Lemma~\ref{lem-cov-concavity}. Let
\[
 K_t=(1-t)K+tL,
 \qquad
 K_t(s)=\{z:g_{K_t}(z)^{1/n}\geq s\}.
\]
Then, for $s\in(0,|K|^{1/n})$ and $r\in(0,|L|^{1/n})$,
\begin{equation}\label{eq:convolution-body-inclusion}
 (1-t)K(s)+tL(r)\subset K_t((1-t)s+tr).
\end{equation}
Indeed, if $z\in K(s)$ and $w\in L(r)$, then $g_K(z)^{1/n}\geq s$ and $ g_L(w)^{1/n}\geq r$. Lemma~\ref{lem-cov-concavity} yields
\[
 g_{K_t}((1-t)z+tw)^{1/n} \geq(1-t)s+tr,
\]
which proves \eqref{eq:convolution-body-inclusion}.

We now prove the transfer principle from Theorem~\ref{thm-intro-transfer}.
% The preceding observations lead to the transfer principle stated in Theorem~\ref{thm-intro-transfer}.

% The preceding observations lead to the following transfer principle from symmetric dual quermassintegrals to chord integrals.
% \begin{thm}%[Transfer from symmetric dual volumes to chord integrals]
% \label{thm:dual-to-chord}
% Let $n\geq2$, $q>1$, and set $\beta=q-1$. Assume that
% \begin{equation}\label{eq:symmetric-dual-input}
%  \widetilde V_\beta((1-t)C+tD)^{\frac{1}{\beta}}
%  \geq (1-t)\widetilde V_\beta(C)^{\frac{1}{\beta}}+t\widetilde V_\beta(D)^{\frac{1}{\beta}}
% \end{equation}
% holds for all $C,D\in\mathcal K_e^n$ and $0\leq t\leq1$.
% Then, for all $K,L\in\mathcal K^n$,
% \begin{equation}\label{eq:transfer-chord-BM}
%  I_q((1-t)K+tL)^{\frac{1}{n+q-1}}
%  \geq (1-t)I_q(K)^{\frac{1}{n+q-1}}  +tI_q(L)^{\frac{1}{n+q-1}}.
% \end{equation}
% % where $K_t=(1-t)K+tL$.
% For $0<t<1$, equality in \eqref{eq:transfer-chord-BM} holds if and only if $K$ and $L$ are homothetic. 
% % In particular, if
% % \eqref{eq:symmetric-dual-input} holds for
% % $0<\beta\leq\beta_0$, then \eqref{eq:transfer-chord-BM} holds for
% % \[
% %  1<q\leq1+\beta_0.
% % \]
% \end{thm}

\begin{proof}[Proof of Theorem~\ref{thm-intro-transfer}]
For $s>0$, let $K_t=(1-t)K+tL$, and recall the convolution bodies
\[
 K(s)=\{z:g_K(z)^{1/n}\geq s\},
 \qquad
 L(r)=\{z:g_L(z)^{1/n}\geq r\}.
\]
They belong to $\mathcal K_e^n$. By \eqref{eq:mu-dual-volume}, assumption \eqref{BM-tildeV} is equivalent to
\[
 \mu_\beta((1-t)C+tD)^{\frac{1}{\beta}}
 \geq (1-t)\mu_\beta(C)^{\frac{1}{\beta}}+t\mu_\beta(D)^{\frac{1}{\beta}}.
\]
% for $C,D\in\mathcal K_e^n$.
It follows from \eqref{eq:convolution-body-inclusion} that
\begin{equation}\label{eq-level-mu}
     \mu_\beta(K_t((1-t)s+tr))^{\frac{1}{\beta}}
 \geq (1-t)\mu_\beta(K(s))^{\frac{1}{\beta}}+t\mu_\beta(L(r))^{\frac{1}{\beta}}.
\end{equation}

Since $g_K(z)\leq g_K(0)=|K|$ and
\[
 g_K(z)
 =n\int_0^{g_K(z)^{\frac{1}{n}}}s^{n-1}\,\dd s
 = n\int_0^\infty s^{n-1}\chi_{K(s)}(z)\,\dd s,
\]
then we have
\[
 \begin{aligned}
 \int_{\mathbb R^n}g_K(z)|z|^{\beta-n}\,\dd z
 &= n\int_0^\infty s^{n-1} \int_{K(s)}|z|^{\beta-n}\,\dd z\,\dd s
 = n\int_0^{|K|^{\frac{1}{n}}} s^{n-1}\mu_\beta(K(s))\,\dd s.
 \end{aligned}
\]
% Here the last upper limit follows from $g_K(z)\leq g_K(0)=|K|$.
Set $f_t(s)=n s^{n-1}\mu_\beta(K_t(s))$ and $p=1/(n+\beta-1)$, extending the height densities by zero outside their positive intervals.
The level-set inequality \eqref{eq-level-mu}  and H\"older's inequality give
\begin{equation*}%\label{eq:level-height-BBL}
 f_t((1-t)s+tr) \geq\big((1-t)f_0(s)^p+t f_1(r)^p\big)^{\frac{1}{p}}
\end{equation*}
whenever $f_0(s)f_1(r)>0$. Applying the one-dimensional Borell-Brascamp-Lieb inequality \cite{BL}, with $p/(1+p)=1/(n+\beta)$, we therefore obtain
\begin{align*}
    \left(\int_{\mathbb R^n}g_{K_t}(z)|z|^{\beta-n}\,\dd z\right)^{\frac{1}{n+\beta} }
    \geq (1-t)\left(\int_{\mathbb R^n}g_{K}(z)|z|^{\beta-n}\,\dd z\right)^{\frac{1}{n+\beta} }+t\left(\int_{\mathbb R^n}g_{L}(z)|z|^{\beta-n}\,\dd z\right)^{\frac{1}{n+\beta} }.
\end{align*}
By $\beta=q-1$, the covariogram representation \eqref{chord-energy} gives \eqref{BM-Iq}.

It remains to prove the equality statement. Let $\psi_t$ be the
jointly concave envelope of $g_K^{1/n}$ and $g_L^{1/n}$
constructed in Section~\ref{sec:weighted-marginals}. The preceding
level-set argument applied to $\psi_t$ gives
\begin{align*}
    \left(\int_{\mathbb R^n}\psi_t(z)^n\,\dd\mu_\beta(z)\right)^{\frac{1}{n+\beta} }
    \geq (1-t)\left(\int_{\mathbb R^n}g_{K}(z)\,\dd\mu_\beta(z)\right)^{\frac{1}{n+\beta} }+t\left(\int_{\mathbb R^n}g_{L}(z)\,\dd\mu_\beta(z)\right)^{\frac{1}{n+\beta} }
\end{align*}
If $K$ and $L$ are not homothetic, the classical Brunn-Minkowski inequality gives
\[
 g_{K_t}(0)^{\frac{1}{n}} = |K_t|^{\frac{1}{n}}
 > (1-t)|K|^{\frac{1}{n}}+t|L|^{\frac{1}{n}} =\psi_t(0).
\]
By continuity, this strict inequality persists on a neighborhood of the origin. Since $\beta>0$, $\mu_\beta$ is positive on every nonempty open set, and therefore
\[
 \int_{\mathbb R^n} \bigl(g_{K_t}(z)-\psi_t(z)^n\bigr)\,\dd\mu_\beta(z)>0.
\]
Thus equality in \eqref{BM-Iq} is impossible unless $K$ and $L$ are homothetic. The converse follows from translation invariance and homogeneity.
\end{proof}

We briefly recall the relevant developments for dual quermassintegrals before applying the transfer principle. 
Xi and Zhang \cite{XZ} proved the Brunn-Minkowski inequality for $\widetilde V_\beta$ in the range $0<\beta\leq1$.
Sadovsky and Zhang \cite{SZ25} established the symmetric Brunn-Minkowski inequality for $\widetilde V_\beta$ in the range $0<\beta\leq n$.
For unconditional bodies, Hu and Ivaki \cite{HI26} further proved the inequality for $n<\beta\leq n+3-2\sqrt{2}$ and at $\beta=n+2$, along with an $L_0$ version for $0<\beta<n$ and $L_p$ extensions for $0<p<1$. 
In our previous work \cite{LWZ}, we extended the range to $n<\beta\leq n+1$ for unconditional bodies and established the endpoint $\beta=n+2$ for origin-symmetric convex bodies. The full range $n<\beta\leq n+1$ for all origin-symmetric bodies then follows from Corollary~\ref{cor-dual-BM}.
% In our previous work \cite{LWZ}, this range was extended to $ n<\beta\leq n+1$ for unconditional bodies, together with the endpoint $\beta=n+2$. For the full class of origin-symmetric bodies, the range $n<\beta\leq n+1$ is provided by Corollary~\ref{cor-dual-BM}.

% \begin{cor}\label{cor:transfer-dual}
% For every $K,L\in\mathcal K^n$ and $q\in [1,n+2]\cup \{n+3\}$,
% \begin{equation}\label{eq:transfer-chord-range}
%  I_q(K_t)^{\frac1{n+q-1}}
%  \geq (1-t)I_q(K)^{\frac1{n+q-1}} +tI_q(L)^{\frac1{n+q-1}},
%  \qquad 0<t<1.
% \end{equation}
% Equality holds if and only if $K$ and $L$ are homothetic.
% \end{cor}

Since $q=\beta+1$, Theorem~\ref{thm-intro-transfer} yields the chord Brunn-Minkowski inequality for $1<q\le n+2$ and at $q=n+3$. Together with the classical case $q=1$, this gives an alternative proof of the general-dimensional part of Theorem~\ref{thm-main}.
In dimension two, Tao's result~\cite{Tao26} supplies the remaining range $4<q<5$ through Theorem~\ref{thm-intro-transfer}, completing the proof of Theorem~\ref{thm-main}.

% \begin{cor}\label{cor:planar-five}
% For planar convex bodies $K,L$ and $q\in[1,5]$,
% \begin{equation}\label{eq:planar-five}
%  I_q(K_t)^{\frac1{q+1}}
%  \geq (1-t)I_q(K)^{\frac1{q+1}} +tI_q(L)^{\frac1{q+1}},
%  \qquad 0<t<1.
% \end{equation}
% Equality holds if and only if $K$ and $L$ are homothetic.
% \end{cor}

\section{Proof of Theorem \ref{thm:fractional-main}}
\label{sec:fractional-simplices}

In this section, we focus on the range $0<q<1$.
A useful observation of Cheung and Xiong \cite{CX} is that, for a simplex $S\in\mathcal{K}^n$,
\begin{equation}\label{eq:simplex-covariogram}
 g_S(ru)
 = |S|\left(1-\frac{r}{\rho_{S-S}(u)}\right)_+^n,
 \qquad r\geq0.
\end{equation}
Thus $g_S^{1/n}$ is affine along each ray until it vanishes.
This suggests comparing a general covariogram with its tangent at the origin. The sharp comparison with $T_q$, combined with concavity of $T_q^{1/(n+q-1)}$, yields the fractional inequality for simplices.

\subsection{The sharp comparison between $I_q$ and $T_q$}\,

For $u\in\mathbb S^{n-1}$, write
\[
 J_q(K,u)=\int_{u^\perp}X_K(y,u)^q\,\dd y,
 \qquad
 I_q(K)=\int_{\mathbb S^{n-1}}J_q(K,u)\,\dd\sigma(u).
\]
Slicing $K$ by lines parallel to $u$ gives
\begin{equation}\label{eq:frac-covariogram-slicing}
 g_K(ru) = \int_{u^\perp}(X_K(y,u)-r)_+\,\dd y,
 \qquad r\geq0.
\end{equation}
Consequently,
\[
 |K|-g_K(ru)= \int_{u^\perp}\min\{X_K(y,u),r\}\,\dd y.
\]
Since $a^q = q(1-q)\int_0^\infty \min\{a,r\}r^{q-2}\,\dd r$ for $0<q<1$ and $a\geq 0$, we obtain
\begin{equation}\label{eq:frac-directional-deficit}
 J_q(K,u)
 = q(1-q)\int_0^\infty \bigl(|K|-g_K(ru)\bigr)r^{q-2}\,\dd r.
\end{equation}

Recall that $g_K^{1/n}$ is concave on $K-K$. Moreover, it follows from \eqref{eq:frac-covariogram-slicing} that
\[
 \left.\frac{\dd}{\dd r}\right|_{r=0+}g_K(ru)
 =-\int_{u^\perp}\chi_{\{X_K(y,u)>0\}}\,\dd y
 =-b_K(u).
\]
Thus $g_K(ru)^{1/n}$ lies below its tangent at the origin, which gives
\begin{equation}\label{eq:frac-tangent-bound}
 g_K(ru)
 \leq |K| \left(1-\frac{b_K(u)}{n|K|}r\right)_+^n.
\end{equation}
For a simplex $S$, using \eqref{eq:simplex-covariogram} and $\rho_{S-S}(u)=\frac{n|S|}{b_S(u)}$, we obtain equality in \eqref{eq:frac-tangent-bound}.

We next obtain the sharp comparison with the functional $T_q$. Recall that
\[
 T_q(K)
 = |K|^q\int_{\mathbb S^{n-1}} b_K(u)^{1-q}\,\dd\sigma(u),
 \qquad
 c_{n,q}
 = \frac{n^q\Gamma(q+1)\Gamma(n)}{\Gamma(n+q)}.
\]

\begin{prop}\label{prop:frac-directional-comparison}
For every $K\in\mathcal K^n$ and
$u\in\mathbb S^{n-1}$,
\begin{equation}\label{eq:frac-directional-comparison}
 c_{n,q}|K|^q b_K(u)^{1-q}
 \leq J_q(K,u)
 \leq |K|^q b_K(u)^{1-q}.
\end{equation}
The lower bound is attained in every direction by every simplex.
\end{prop}

\begin{proof}
The upper bound follows from Jensen's inequality:
\[
 \begin{aligned}
 J_q(K,u)
 &\leq b_K(u) \left(\frac1{b_K(u)}\int_{u^\perp}X_K(y,u)\,\dd y\right)^q
 =|K|^q b_K(u)^{1-q}.
 \end{aligned}
\]

For the lower bound, set $ R=\frac{n|K|}{b_K(u)}.$
Combining \eqref{eq:frac-directional-deficit} with \eqref{eq:frac-tangent-bound} gives
\[
 \begin{aligned}
 J_q(K,u)
 &\geq q(1-q)|K| \int_0^\infty \left[1-\left(1-\frac rR\right)_+^n\right] r^{q-2}\,\dd r\\
 &= q(1-q)|K|R^{q-1} \int_0^\infty \bigl[1-(1-s)_+^n\bigr]s^{q-2}\,\dd s\\
 &= q n B(q,n) |K|R^{q-1}
 =c_{n,q}|K|^q b_K(u)^{1-q}.
 \end{aligned}
\]
The equality statement follows from the equality case of \eqref{eq:frac-tangent-bound}.
\end{proof}

Integrating \eqref{eq:frac-directional-comparison} over $\mathbb S^{n-1}$ gives
\begin{equation*}%\label{eq:frac-sharp-comparison}
 c_{n,q}T_q(K)
 \leq I_q(K)
 \leq T_q(K).
\end{equation*}
The lower constant is sharp, since every simplex attains equality.
Moreover, $0<c_{n,q}<1.$

Conversely, suppose that $I_q(K)=c_{n,q}T_q(K)$. Then the lower
bound in \eqref{eq:frac-directional-comparison} is attained for
almost every $u$. By \eqref{eq:frac-directional-deficit}, the
nonnegative difference in \eqref{eq:frac-tangent-bound} has zero
integral against the positive weight $r^{q-2}$. Hence, by continuity
in $r$ and $u$,
\[
 g_K(ru)=|K|\left(1-\frac{r b_K(u)}{n|K|}\right)_+^n,
 \qquad r\geq0,\quad u\in\mathbb S^{n-1}.
\]
Let $D_K=K-K$ be the difference body of $K$. Since $\operatorname{supp}g_K=D_K$, the comparison of radial supports gives $\rho_{D_K}(u)=n|K|/b_K(u)$. Thus, we obtain
\[
 |K|^2
 =\int_{D_K} g_K(z)\,\dd z
 =n|K||D_K|\int_0^1(1-s)^n s^{n-1}\,\dd s
 =\frac{|K||D_K|}{\binom{2n}{n}}.
\]
By the equality case of the Rogers-Shephard difference-body inequality \cite[Theorem 10.1.4]{Schneider}, $K$ is a simplex. This completes the proof of Theorem~\ref{thm:fractional-main} (1).

\subsection{Concavity of the comparison functional $T_q$}\,

In the range $0<q<1$, our comparison argument does not yield the concavity of $I_q$ itself. We therefore study instead the comparison functional $T_q$. The key step is to establish the Brunn-Minkowski inequality for
\[
 A_q(K)
 = \left(\int_{\mathbb S^{n-1}} b_K(u)^{1-q}\,\dd\sigma(u) \right)^{\frac{1}{1-q}}.
\]
Since $b_K(u)=|P_{u^\perp}K|_{n-1}$, Cauchy's projection formula allows us to represent $A_q$ as an extremal mixed-volume functional.
More precisely, we have the following representation.

\begin{prop}\label{prop:Aq-mixed-volume}
For every $K\in\mathcal K^n$ and $0<q<1$,
\begin{equation}\label{eq:Aq-mixed-volume}
 A_q(K)
 =\min\left\{ nV_1(K,Z_w):\,w\in C(\mathbb{S}^{n-1}),\, w>0, \, \int_{\mathbb S^{n-1}}w(u)^{\frac{q-1}{q}}\,\dd\sigma(u)=1\right\}
 % \min_{\substack{w>0\\
 % \int_{\mathbb S^{n-1}}
 % w(u)^{(q-1)/q}\,\dd\sigma(u)=1}}
 % nV(K[n-1],Z_w),
\end{equation}
where $Z_w$ is the centrally symmetric zonoid defined by
\begin{equation}\label{eq:Zw-support}
 h_{Z_w}(v)
 =\frac12\int_{\mathbb S^{n-1}} w(u)|\langle u,v\rangle|\,\dd\sigma(u).
\end{equation}
\end{prop}

\begin{proof}
By the weighted power-mean identity, we have
\[
 A_q(K)=\min\left\{\int_{\mathbb S^{n-1}}w(u)b_K(u)\,\dd\sigma(u):\, w\in C(\mathbb{S}^{n-1}),\, w>0, \, \int_{\mathbb S^{n-1}}w(u)^{\frac{q-1}{q}}\,\dd\sigma(u)=1\right\}.
\]
The minimum is attained at
\[
 w_K(u)
 = \left( \int_{\mathbb S^{n-1}} b_K(v)^{1-q}\,\dd\sigma(v) \right)^{\frac{q}{1-q}} b_K(u)^{-q}.
\]
By Cauchy's projection formula and the definition of $Z_w$, we obtain
\[
 \int_{\mathbb S^{n-1}}w(u)b_K(u)\,\dd\sigma(u)
 = nV_1(K,Z_w).
\]
Substitution into the power-mean identity proves \eqref{eq:Aq-mixed-volume}.
\end{proof}

\begin{prop}\label{prop:Aq-BM}
For every $0<q<1$ and $K,L\in\mathcal K^n$,
\begin{equation}\label{eq:Aq-BM}
 A_q((1-t)K+tL)^{\frac{1}{n-1}}
 \geq (1-t)A_q(K)^{\frac{1}{n-1}} +tA_q(L)^{\frac{1}{n-1}},
 \qquad 0<t<1.
\end{equation}
Consequently, we have
\begin{equation}\label{eq:frac-T-BM}
 T_q((1-t)K+tL)^{\frac{1}{n+q-1}}
 \geq (1-t)T_q(K)^{\frac{1}{n+q-1}}+tT_q(L)^{\frac{1}{n+q-1}},
  \qquad 0<t<1.
\end{equation}
Equality holds in \eqref{eq:frac-T-BM} if and only if $K$ and $L$ are homothetic.
\end{prop}

\begin{proof}
For an admissible weight $w$, the Aleksandrov-Fenchel inequality gives
\[
V_1((1-t)K+tL,Z_w)^{\frac{1}{n-1}}\geq (1-t) V_1(K,Z_w)^{\frac{1}{n-1}}+t V_1(L,Z_w)^{\frac{1}{n-1}}
\]
Choose a minimizing weight for $(1-t)K+tL$. Since the same weight is admissible for $K$ and $L$, \eqref{eq:Aq-mixed-volume} gives \eqref{eq:Aq-BM}.
Furthermore, we have
\[
 T_q(K)^{\frac{1}{n+q-1}}
 =\bigl(|K|^{\frac{1}{n}}\bigr)^{\frac{nq}{n+q-1}}
 \bigl(A_q(K)^{\frac{1}{n-1}}\bigr)^{ \frac{(n-1)(1-q)}{n+q-1}}.
\]
Thus, the classical Brunn-Minkowski inequality, \eqref{eq:Aq-BM}, and Proposition~\ref{prop-calculus}(iii) yield \eqref{eq:frac-T-BM}.

Suppose equality holds in \eqref{eq:frac-T-BM}. Since the weighted geometric mean is strictly increasing in each variable, equality must hold in the classical Brunn-Minkowski inequality. Hence $K$ and $L$ are homothetic. The converse is immediate.
\end{proof}

Finally, Proposition~\ref{prop:frac-directional-comparison} and Proposition~\ref{prop:Aq-BM} together prove Theorem~\ref{thm:fractional-main}.

\section{Minkowski inequalities and uniqueness}
\label{sec:consequences}

In this section, we derive Minkowski inequalities and use them to prove the uniqueness results for chord measures. Throughout this section, set $d=n+q-1$.

\subsection{The case \texorpdfstring{$p=1$}{p=1}}\,

Recall the first-variation formula in \cite{LXYZ}
\begin{equation}\label{chord-variation}
 \left.\frac{\dd}{\dd s}\right|_{s=0^+}I_q(K+sL)
 =\int_{\mathbb S^{n-1}}h_L(u)\,\dd F_q(K,u).
\end{equation}
Define the mixed chord integral by
\[
 I_q(K,L)
 =\frac1d\int_{\mathbb S^{n-1}}h_L(u)\,\dd F_q(K,u).
\]
By the homogeneity of $I_q$, we have
\begin{equation}\label{Euler}
 I_q(K,K)=I_q(K).
\end{equation}

The Brunn-Minkowski inequality for $I_q$ now yields the following Minkowski-type inequality, as well as the uniqueness of the chord measure.

\begin{thm}\label{thm-mixed-Min}
Let $K,L\in\mathcal K^n$ and
$q\in[1,n+2]\cup\{n+3\}$. Then
\begin{equation}\label{chord-Minkowski}
 I_q(K,L)
 \geq I_q(K)^{\frac{d-1}{d}}I_q(L)^{\frac1d}.
\end{equation}
Equality holds if and only if $K$ and $L$ are homothetic.

In the planar case, the same assertions hold for every $1\leq q\leq5$.
\end{thm}

\begin{proof}
For $s\geq0$, homogeneity and the chord Brunn-Minkowski inequality give
\[
 I_q(K+sL)^{\frac{1}{d}}
 \geq I_q(K)^{\frac{1}{d}}+sI_q(L)^{\frac{1}{d}}.
\]
Differentiating at $s=0+$ and using
\eqref{chord-variation} and \eqref{Euler}, we obtain
\[
 I_q(K,L)
 = \frac1d\left.\frac{\dd}{\dd s}\right|_{s=0+} I_q(K+sL)
 \geq I_q(K)^{\frac{d-1}{d}}I_q(L)^{\frac1d}.
\]
If equality holds, set $f(s)=I_q(K+sL)^{1/d}$. Then
$f'_+(0)=I_q(L)^{1/d}$. By concavity and the preceding inequality,
\[
 f(0)+sI_q(L)^{1/d}\leq f(s)\leq f(0)+sf'_+(0),
 \qquad s\geq0.
\]
Hence $f$ is affine. Then the equality characterization in the chord Brunn-Minkowski inequality implies that $K$ and $L$ are homothetic. The converse follows from homogeneity and translation invariance.

The planar assertion follows from the planar part of Theorem \ref{thm-main}.
\end{proof}

\subsection{The case \texorpdfstring{$p>1$}{p>1}}\,

Let $K,L\in\mathcal K_o^n$. Recall that Firey's $p$-sum is defined by
\[
 h_{(1-t)\cdot_pK+_pt\cdot_pL}^p
 =(1-t)h_K^p+th_L^p.
\]
Since $I_q$ is strictly increasing under proper containment, the chord Brunn-Minkowski inequality and Proposition~\ref{prop-calculus}(iv) give the following.

\begin{thm}\label{thm-Firey}
Let $p>1$ and $q\in[1,n+2]\cup\{n+3\}$. For $K,L\in\mathcal K_o^n$ and $0< t<1$,
\begin{equation}\label{Lp-chord-BM}
 I_q((1-t)\cdot_pK+_pt\cdot_pL)^{\frac{p}{d}}
 \geq (1-t)I_q(K)^{\frac{p}{d}}+tI_q(L)^{\frac{p}{d}}.
\end{equation}
Equality holds if and only if $K$ and $L$ are dilates.

In the planar case, the same assertions hold for every $1\leq q\leq5$.
\end{thm}

For the $L_p$ chord measure
\[
 \dd F_{p,q}(K,u)
 = h_K(u)^{1-p}\,\dd F_q(K,u),
\]
define the $L_p$ mixed chord integral by
\[
 I_{p,q}(K,L)
 = \frac1d\int_{\mathbb S^{n-1}}h_L(u)^p\, \dd F_{p,q}(K,u).
\]
Then $I_{p,q}(K,K)=I_q(K)$, and the first-variation formula
\cite{LXYZ} gives
\begin{equation}\label{Lp-chord-variation}
 \left.\frac{\dd}{\dd s}\right|_{s=0^+} I_q(K+_ps\cdot_pL)
 = \frac1p\int_{\mathbb S^{n-1}}h_L(u)^p\,
 \dd F_{p,q}(K,u).
\end{equation}
Differentiating the $L_p$ chord Brunn-Minkowski inequality gives the $L_p$ Minkowski inequality.

\begin{thm}\label{thm-Lp-mixed}
Let $K,L\in\mathcal K_o^n$, $p>1$, and $ q\in[1,n+2]\cup\{n+3\}.$
Then
\begin{equation}\label{Lp-mixed-inequality}
 I_{p,q}(K,L)
 \geq I_q(K)^{\frac{d-p}{d}} I_q(L)^{\frac{p}{d}}.
\end{equation}
Equality holds if and only if $K$ and $L$ are dilates.

In the planar case, the same assertions hold for every $1\leq q\leq 5$.
\end{thm}

Indeed, \eqref{Lp-mixed-inequality} and its equality characterization follow by differentiating \eqref{Lp-chord-BM} at $s=0^+$ and using \eqref{Lp-chord-variation}, exactly as in the case $p=1$.

\subsection{Proof of Theorem~\ref{thm-chord-uniqueness}}\,

We now prove the uniqueness statements for chord measures. They follow by applying the two Minkowski inequalities above in opposite directions.

\begin{proof}[Proof of Theorem~\ref{thm-chord-uniqueness}]
We first consider the case $p=1$. Suppose that 
$$F_q(K,\cdot)=F_q(L,\cdot).$$
By the definition of the mixed chord integral and \eqref{Euler}, we have
\[
I_q(K,L)=I_q(L),
\qquad
I_q(L,K)=I_q(K).
\]
Applying \eqref{chord-Minkowski} in both directions gives
\[
I_q(L)\geq I_q(K)^{\frac{d-1}{d}}I_q(L)^{\frac1d},
\qquad
I_q(K)\geq I_q(L)^{\frac{d-1}{d}}I_q(K)^{\frac1d}.
\]
Since $d=n+q-1>1$, these inequalities imply $ I_q(K)\leq I_q(L)$ and $I_q(L)\leq I_q(K)$, then $I_q(K)=I_q(L)$. 
Equality holds in \eqref{chord-Minkowski}, so $K$ and $L$ are homothetic. Using $I_q(K)=I_q(L)$ again, we conclude that $K$ and $L$ differ only by a translation.

Now let $p>1$ and suppose that
$$F_{p,q}(K,\cdot)=F_{p,q}(L,\cdot).$$
Similarly, we get
\[
I_{p,q}(K,L)=I_q(L),
\qquad
I_{p,q}(L,K)=I_q(K).
\]
Applying \eqref{Lp-mixed-inequality} in both directions yields
\[
I_q(L)^{\frac{d-p}{d}} \geq I_q(K)^{\frac{d-p}{d}},
\qquad
I_q(K)^{\frac{d-p}{d}} \geq I_q(L)^{\frac{d-p}{d}}.
\]

If $p\neq d$, these inequalities imply $I_q(K)=I_q(L)$. Equality in \eqref{Lp-mixed-inequality} then shows that $K$ and $L$ are dilates, and $I_q(K)=I_q(L)$ forces $K=L$.
If $p=d$, equality holds in \eqref{Lp-mixed-inequality}, so $K$ and $L$ are dilates.

The planar assertions follow from the planar chord Brunn-Minkowski range.
\end{proof}

\section{Extensions to affine sections and mixed energies}\label{sec:extensions}

The preceding arguments extend in two directions: from line sections to higher-dimensional affine sections, and from self-covariograms to mixed covariograms.

\subsection{Affine-section integrals}\,

Huang, Yang, and Zhang \cite[Section~7.6]{HYZ} suggest extending the study of chord integrals to intersections with higher-dimensional affine planes. Following this perspective, we consider moments of affine section volumes, rather than only over subspaces through the origin.
% Following the perspective in \cite[Section~7.6]{HYZ}, we consider
% moments of volumes of affine sections.

For $1\leq m<n$ and $q>0$, set
\[
 a_K(y,E)=|K\cap(y+E)|_m,
 \qquad E\in G(n,m),\quad y\in E^\perp,
\]
and define
\begin{equation}\label{def-Imq}
 I_{m,q}(K)
 =\int_{G(n,m)}\int_{E^\perp}a_K(y,E)^q\,\dd y\,\dd\nu(E), \qquad q>0,
\end{equation}
where $\nu$ is probability Haar measure on $G(n,m)$. At $q=0$, we use the convention
\begin{equation}\label{def-Im0}
 I_{m,0}(K)
 = \int_{G(n,m)}|P_{E^\perp}K|_{n-m}\,\dd\nu(E).
\end{equation}
It follows that $I_{m,q}$ is translation invariant and homogeneous of degree $n-m+mq$.
% Unlike averages of sections through a prescribed origin, %centro-section...
% $I_{m,q}$ is translation invariant.
%Here ``affine'' refers to the planes being integrated, not to invariance under general affine maps.
Moreover,
\begin{equation}\label{Imq-hom}
 % I_{m,q}(rK)=r^{n-m+mq}I_{m,q}(K),
 % \qquad
 I_{1,q}=I_q,\qquad
 I_{m,1}(K)=|K|,\qquad
  I_{m,0}(K)
 =\frac{\omega_{n-m}\omega_m}{\binom{n}{m}\omega_n}V_{n-m}(K),
\end{equation}
where the last identity is Kubota's formula \cite{Schneider}, and the intrinsic volumes are normalized by
\[
 |K+rB^n| = \sum_{j=0}^n\omega_{n-j}V_j(K)r^{n-j}.
\]
At the second moment, Zhang's formula \eqref{intro-Zhang} yields
\begin{equation}\label{Im2-chord}
 I_{m,2}(K)
 = \frac{m\omega_m}{n\omega_n} \mathcal E_{m-n}(K)
 =\frac{\omega_m}{m+1}I_{m+1}(K).
\end{equation}

\begin{thm}\label{thm-affine}
Let $1\leq m<n$ and $q\in[1,1+\frac1m]\cup\{0,2\}$. For $K,L\in\mathcal K^n$,
\begin{equation}\label{BM-Imq}
 I_{m,q}((1-t)K+tL)^{\frac{1}{n-m+mq}}
 \geq (1-t)I_{m,q}(K)^{\frac{1}{n-m+mq}} +tI_{m,q}(L)^{\frac{1}{n-m+mq}},\,\, 0<t<1.
\end{equation}
Equality holds if and only if $K$ and $L$ are homothetic, unless $m=n-1$ and $q=0$, in which case the inequality is an identity.
\end{thm}

\begin{proof}
For $1<q\leq1+1/m$ and $0<t<1$, set $\alpha=m(q-1)\in(0,1]$ and $K_t=(1-t)K+tL$. 
We use the convex-lifting argument from the previous section, combining a lift to one higher dimension with the linear Brunn-Minkowski inequality.

Fix $E\in G(n,m)$. For $y\in E^\perp$, define the section-root profile
\[
 r_{K,E}(y)
 =a_K(y,E)^{\frac{1}{m}}
 =|K\cap(y+E)|_m^{\frac{1}{m}}.
\]
Brunn's section theorem implies that $r_{K,E}$ is concave on $P_{E^\perp}K$. Moreover, section inclusion and the $m$-dimensional Brunn-Minkowski inequality give
\[
 r_{K_t,E}((1-t)y+tz)
 \geq (1-t)r_{K,E}(y)+t r_{L,E}(z).
\]
Since $0<\alpha\leq1$, the function $s\mapsto s^\alpha$ is increasing and concave. Hence
\[
 r_{K_t,E}((1-t)y+tz)^\alpha
 \geq \bigl((1-t)r_{K,E}(y)+t r_{L,E}(z)\bigr)^\alpha
 \geq (1-t)r_{K,E}(y)^\alpha+t r_{L,E}(z)^\alpha.
\]
Consequently, the hypograph
\[
 \mathcal L_{E,q}(K)
 = \left\{ (x,s)\in\mathbb R^n\times\mathbb R: x\in K,\quad 0\leq s\leq r_{K,E}(P_{E^\perp}x)^\alpha \right\}
\]
is a convex body in $\mathbb R^{n+1}$. The preceding inequality also gives the lift inclusion
\[
 (1-t)\mathcal L_{E,q}(K) +t\mathcal L_{E,q}(L)
 \subset \mathcal L_{E,q}(K_t).
\]
Moreover, we have
\[
 \begin{aligned}
 |\mathcal L_{E,q}(K)|
 &=\int_K r_{K,E}(P_{E^\perp}x)^\alpha\,\dd x\\
 &=\int_{E^\perp}r_{K,E}(y)^\alpha
   |K\cap(y+E)|_m\,\dd y\\
 &=\int_{E^\perp}a_K(y,E)^q\,\dd y.
 \end{aligned}
\]

Suppose first that $|K|=|L|$. The projections of $\mathcal L_{E,q}(K)$ and $\mathcal L_{E,q}(L)$ onto $\mathbb R^n$ are $K$ and $L$, respectively. Thus the linear Brunn-Minkowski inequality in Theorem~\ref{thm-Bonnesen}, together with the lift inclusion, gives
\[
 \int_{E^\perp}a_{K_t}(y,E)^q\,\dd y
 \geq(1-t)\int_{E^\perp}a_K(y,E)^q\,\dd y +t\int_{E^\perp}a_L(y,E)^q\,\dd y.
\]
Integrating over $E\in G(n,m)$ yields the equal-volume linear inequality
\[
 I_{m,q}(K_t)
 \geq (1-t)I_{m,q}(K)+tI_{m,q}(L).
\]
Since $I_{m,q}$ is homogeneous of degree $n-m+mq$, Proposition~\ref{prop-calculus}(ii) gives \eqref{BM-Imq} for all $K,L$.

For the equality case, normalize $K$ and $L$ to have equal volume. Equality in the averaged inequality forces equality in the directional inequalities for almost every $E$. The equality statement in Theorem~\ref{thm-Bonnesen} then shows that the normalized projections $K$ and $L$ are translates. Without the normalization, $K$ and $L$ must be homothetic. 
The converse follows from homogeneity and translation invariance.

Finally, the case $q=1$ is the classical volume Brunn-Minkowski inequality.
For $q=0$, \eqref{Imq-hom} reduces the assertion to $V_{n-m}^{1/(n-m)}$. If $n-m\geq2$, the inequality and equality case follow from the Aleksandrov-Fenchel inequalities.
If $n-m=1$, $V_1$ is Minkowski linear and every pair gives equality.
For $q=2$, \eqref{Im2-chord} reduces the assertion to the chord Brunn-Minkowski inequality at $m+1$. Since $2\leq m+1\leq n$, Theorem~\ref{thm-main} gives both the inequality and its equality
characterization.
\end{proof}

\subsection{Mixed covariograms and energies}\,

The covariogram argument extends naturally to two bodies; see, e.g., \cite{AG}.
For $K,L\in\mathcal K^n$, define the mixed covariogram
\begin{equation}\label{mixed-g}
 g_{K,L}(z)
 = |K\cap(L+z)|
 = (\chi_K*\chi_{-L})(z).
\end{equation}
Then $g_{K,K}=g_K$, $\operatorname{supp}g_{K,L}=K-L$, and
\begin{equation}\label{mixed-energy}
 \mathcal E_\beta(K,L)
 =\int_{\mathbb R^n}g_{K,L}(z)|z|^\beta\,\dd z
 =\int_K\int_L|x-y|^\beta\,\dd x\,\dd y,
 \qquad \beta>-n.
\end{equation}
For two pairs $(K_i,L_i)$, $i=0,1$, write
\[
 K_t=(1-t)K_0+tK_1,
 \qquad
 L_t=(1-t)L_0+tL_1.
\]
We call the pairs $(K_0,L_0)$ and $(K_1,L_1)$ simultaneously homothetic if there exist $r>0$ and $a\in\mathbb R^n$ such that
\[
 K_1=rK_0+a,\qquad L_1=rL_0+a.
\]

\begin{lem}\label{lem-mixed-joint}
The domain
\[
 \Omega=\{(t,z):0\leq t\leq1,\ z\in K_t-L_t\}
\]
is convex, and $(t,z)\mapsto g_{K_t,L_t}(z)^{1/n}$ is concave on $\Omega$.
\end{lem}

\begin{proof}
For $z_i\in K_i-L_i$, section inclusion gives
\begin{align*}
    (1-t)\bigl(K_0\cap(L_0+z_0)\bigr) +t\bigl(K_1\cap(L_1+z_1)\bigr)
    \subset  K_t\cap\bigl(L_t+(1-t)z_0+tz_1\bigr).
\end{align*}
% \[
%  \begin{split}
%  &(1-t)\bigl(K_0\cap(L_0+z_0)\bigr)
%  +t\bigl(K_1\cap(L_1+z_1)\bigr)\\
%  &\qquad\subset
%  K_t\cap\bigl(L_t+(1-t)z_0+tz_1\bigr).
%  \end{split}
% \]
The same inclusion between any two intermediate parameters proves convexity of $\Omega$ and joint concavity by the classical Brunn-Minkowski inequality.
\end{proof}

Unlike $g_K$, a mixed covariogram need not be even. Evenness holds when both bodies are origin-symmetric, or on a diagonal path $L_t=K_t$.
It also forces the component centroids to agree, since
\begin{equation}\label{eq:mixed-first-moment}
 \int_{\mathbb R^n}z\,g_{K,L}(z)\,\dd z
 =|K||L|(c_K-c_L),\qquad
 \int_{\mathbb R^n}g_{K,L}(z)\,\dd z=|K||L|.
\end{equation}
Under the evenness assumption, the weighted-marginal and symmetric-lift arguments extend as follows.

% Unlike $g_K$, a mixed covariogram need not be even. Evenness holds when both bodies are origin-symmetric, or on a diagonal path $L_t=K_t$. Under this condition, the weighted-marginal and symmetric-lift arguments extend as follows.

\begin{thm}\label{thm-mixed}
Suppose $g_{K_t,L_t}$ is even for every $0\leq t\leq1$.
For $\beta\in(-n,1]\cup\{2\}$,
\begin{equation}\label{mixed-BM}
 \mathcal E_\beta(K_t,L_t)^{\frac{1}{2n+\beta}}
 \geq (1-t)\mathcal E_\beta(K_0,L_0)^{\frac{1}{2n+\beta}} +t\mathcal E_\beta(K_1,L_1)^{\frac{1}{2n+\beta}}.
\end{equation}
For $\beta\in[0,1]\cup\{2\}$ and $0<t<1$, equality holds if and only if the two pairs are simultaneously homothetic.
%(-n,1-n]\cup
%unless $\beta=0$, in which case independent translations of the two components are allowed.
In particular, this holds when all four component bodies are origin-symmetric.

In the planar case, the inequality and the equality characterization hold for $-2<\beta\leq2$.
\end{thm}

\begin{proof}
For $-n<\beta<0$, we use the regularized weighted-marginal argument from the proof of Theorem~\ref{thm-negative-energy}.
By Lemma~\ref{lem-mixed-joint}, the function $(t,z)\mapsto g_{K_t,L_t}(z)^{1/n}$ is jointly concave. By the evenness of $g_{K_t,L_t}$, applying Theorem~\ref{thm-WM} gives the concavity of
\[
 t\longmapsto
 \left( \int_{\mathbb R^n} g_{K_t,L_t}(z)W_\varepsilon(z)\,\dd z \right)^{\frac{1}{2n}}.
\]
As $\varepsilon\to0^+$, we obtain the concavity of $ t\mapsto\mathcal E_\beta(K_t,L_t)^{1/(2n)}.$
Since $\mathcal E_\beta(K,L)$ is homogeneous of degree $2n+\beta$ under common dilations of $K$ and $L$, Proposition~\ref{prop-calculus}(i) gives \eqref{mixed-BM}.

For $\beta=0$, we have $\mathcal E_0(K,L)=|K||L|$. The classical Brunn-Minkowski inequality for each component and Proposition~\ref{prop-calculus}(iii) give \eqref{mixed-BM}.
Equality forces homothety in each component, while equality in the geometric-mean inequality makes the two dilation factors equal.
Thus equality holds precisely when
\[
 K_1=rK_0+a,\qquad L_1=rL_0+b
\]
for some $r>0$ and $a,b\in\mathbb R^n$. 
By evenness and \eqref{eq:mixed-first-moment}, $c_{K_i}=c_{L_i}$ for $i=0,1$, so $a=b$.
Hence equality holds precisely when the two pairs are simultaneously homothetic.
% The translations are independent because $\mathcal E_0$ depends only on the two volumes.

For $0<\beta\leq1$, define
\[
 \mathcal B(K,L)
 = \left\{ (z,w)\in\mathbb R^n\times\mathbb R^n: z\in K-L,\quad |w|\leq g_{K,L}(z)^{\frac{1}{n}} \right\}.
\]
Since $g_{K,L}$ is even, $\mathcal B(K,L)$ is an origin-symmetric convex body in $\mathbb R^{2n}$. Lemma~\ref{lem-mixed-joint} gives
\[
 (1-t)\mathcal B(K_0,L_0)+t\mathcal B(K_1,L_1)
 \subset \mathcal B(K_t,L_t).
\]
For $P(z,w)=z$, we have
\[
 \int_{\mathcal B(K,L)}|P(z,w)|^\beta\,\dd z\,\dd w
 =\omega_n\mathcal E_\beta(K,L).
\]
Lemma~\ref{lem-projected} therefore proves \eqref{mixed-BM}. 
For the equality case, normalize the endpoint energies to one, with the corresponding change of interpolation parameter.
Projected-moment rigidity and strict monotonicity under proper containment give
\[
 \mathcal B(K_0,L_0) =\mathcal B(K_1,L_1) =\mathcal B(K_t,L_t).
\]
Their equal volumes imply
\[
 |K_0||L_0|=|K_1||L_1|=|K_t||L_t|.
\]
The $\beta=0$ equality case and \eqref{eq:mixed-first-moment} now give simultaneous homothety. Undoing the normalization proves the assertion for the original pairs. The converse is immediate.

% For the equality case, the projected-moment rigidity implies that the two lifts are dilates. The lift inclusion and the fiber description then force the same dilation factor and translation in both components. Hence the pairs are simultaneously homothetic. The converse is immediate.
% The converse follows from homogeneity and simultaneous translation invariance.

For $\beta=2$, we write $\mathcal I_N$ for centered polar inertia in $\mathbb R^N$. Expansion around the two centroids gives
\begin{equation}\label{eq:mixed-quadratic-product}
 \begin{aligned}
 \mathcal E_2(K,L)
 &=|L|\mathcal I_n(K)+|K|\mathcal I_n(L) +|K||L||c_K-c_L|^2,\\
 \mathcal I_{2n}(K\times L)
 &=|L|\mathcal I_n(K)+|K|\mathcal I_n(L).
 \end{aligned}
\end{equation}
By evenness and \eqref{eq:mixed-first-moment}, $c_{K_t}=c_{L_t}$, so
\[
 \mathcal E_2(K_t,L_t) =\mathcal I_{2n}(K_t\times L_t).
\]
Since $K_t\times L_t =(1-t)(K_0\times L_0)+t(K_1\times L_1)$, Hadwiger's inertia inequality in dimension $2n$ gives \eqref{mixed-BM}. Its equality characterization implies $ K_1=rK_0+a$ and $L_1=rL_0+b$ for some $r>0$. The centroid identity gives $a=b$, proving simultaneous homothety.

The level-set transfer argument also applies to the mixed covariograms. More generally, put $\sigma=n+\beta>0$.
By evenness and Lemma~\ref{lem-mixed-joint}, their positive superlevel sets are origin-symmetric convex bodies and satisfy
\[
 (1-t)\{g_{K_0,L_0}^{1/n}\geq s\}
 +t\{g_{K_1,L_1}^{1/n}\geq r\}\subset
 \{g_{K_t,L_t}^{1/n}\geq(1-t)s+tr\}.
\]
Whenever $\mu_\sigma^{1/\sigma}$ satisfies the Brunn-Minkowski inequality on these sets, the integration argument of Section~\ref{sec:dual-to-chord} gives \eqref{mixed-BM}. 

In dimension two, Tao's planar theorem
\cite{Tao26}, with $\sigma=2+\beta$, supplies the remaining
range $1<\beta<2$.
For the planar equality cases, normalize the endpoint energies to one, so equality gives $\mathcal E_\beta(K_t,L_t)=1$.
Set
\[
 f_j(s)=2s\,\mu_{\sigma}
 \bigl(\{g_{K_j,L_j}^{1/2}\geq s\}\bigr),
 \qquad j=0,1,t.
\]
These height densities have integral one, and the level-set
argument gives
\[
 f_t((1-t)s+tr)
 \geq
 \bigl((1-t)f_0(s)^{1/(3+\beta)}
       +t f_1(r)^{1/(3+\beta)}\bigr)^{3+\beta}.
\]
Let $T$ be the increasing transport between the endpoint
densities, so $f_1(T)T'=f_0$ almost everywhere.
Change of variables and H\"older's inequality yield
\[
 \begin{aligned}
 1=\int f_t
 &\geq\int
 \bigl((1-t)f_0(s)^{\frac{1}{3+\beta}} +t f_1(T(s))^{\frac{1}{3+\beta}}\bigr)^{3+\beta}\bigl((1-t)+tT'(s)\bigr)\,\dd s\\
 &\geq\int\bigl((1-t)f_0(s)^{\frac{1}{4+\beta}}   +t(f_1(T(s))T'(s))^{\frac{1}{4+\beta}}\bigr)^{4+\beta}\,\dd s
 =1.
 \end{aligned}
\]
The changes of variables are justified locally on the positive height intervals, where the densities and their reciprocals are bounded.
All inequalities are equalities almost everywhere.
The last H\"older equality and the transport identity give $f_1(T)=f_0$ and $T'=1$.
Equality in the height-factor H\"older step then gives
\[
 \frac{T(s)}s
 =\left(\frac{s f_1(T(s))}{T(s)f_0(s)}\right)^{\frac{1}{2+\beta}},
\]
hence $T(s)=s$.
The rigidity in \cite[Theorem~1.1]{Tao26}, together with equal level-set measures, now forces the three level sets to coincide almost everywhere. Thus
% The unit integrals exclude additional positive levels.
% Layer cake and continuity therefore give
\[
 g_{K_0,L_0}=g_{K_1,L_1}=g_{K_t,L_t}.
\]
Their total masses are equal, so the $\beta=0$ equality case gives simultaneous homothety.
% Undoing normalization proves the assertion.
The converse follows from homogeneity and simultaneous translation invariance.
% Conversely, simultaneous homotheties give equality by
% common homogeneity and simultaneous translation invariance.

\end{proof}

As a consequence, for $1\leq m<n$, the two-point affine-section identity
\[
 \int_{\mathcal A(n,m)} |K\cap F|_m\,|L\cap F|_m\,\dd F
 = \frac{m\omega_m}{n\omega_n} \mathcal E_{m-n}(K,L)
\]
yields a joint Brunn-Minkowski inequality for origin-symmetric component pairs. On the diagonal $K=L$, this recovers the second-moment inequality for $I_{m,2}$.

\section{Further problems}\label{sec:problems}

\begin{figure}[ht]
\centering
\begin{tikzpicture}[x=1cm,y=1cm]
 % Negative indices: shown for reference only
 \draw[densely dashed,line width=.8pt] (0,0)--(1,0);

 % Unresolved ranges
 \draw[densely dashed,line width=.8pt] (1,0)--(2,0);
 \draw[densely dashed,line width=.8pt] (8.1,0)--(9.3,0);
 \draw[densely dashed,line width=.8pt] (9.3,0)--(11.2,0);

 % Known range
 \draw[line width=1.2pt] (2,0)--(8.1,0);

 % Failure range
 \draw[dash dot,line width=1.1pt,->] (11.2,0)--(12.8,0);

 % Interior markers
 \foreach \x in {3.3,6.8}
   \draw (\x,-.07)--(\x,.07);

 % Known endpoints
 \foreach \x in {1,2,8.1,9.3}
   \fill (\x,0) circle (2.2pt);

 % Excluded endpoint and unresolved threshold
 \draw[gray,fill=white] (0,0) circle (2.2pt);
 \draw[fill=white] (11.2,0) circle (2.2pt);

 % Labels
 \node[below=7pt] at (0,0) {$-1$};
 \node[below=7pt] at (1,0) {$0$};
 \node[below=7pt] at (2,0) {$1$};
 \node[below=7pt] at (3.3,0) {$2$};
 \node[below=7pt] at (6.8,0) {$n+1$};
 \node[below=7pt] at (8.1,0) {$n+2$};
 \node[below=7pt] at (9.3,0) {$n+3$};
 \node[below=7pt] at (11.2,0) {$n+1+\beta_*$};
 \node[above=5pt] at (12.7,0) {$q$};

 % Legend
 \draw[line width=1.2pt] (1,-1)--(1.8,-1);
 \node[right] at (1.8,-1) {\small valid};
 \draw[densely dashed,line width=.8pt] (4.2,-1)--(5,-1);
 \node[right] at (5,-1) {\small unresolved here};
 \draw[dash dot,line width=1.1pt] (8.8,-1)--(9.6,-1);
 \node[right] at (9.6,-1) {\small counterexamples};
\end{tikzpicture}
\caption{Parameter ranges for the chord Brunn-Minkowski inequality.\\
For general convex bodies, the inequality is proved on
$[1,n+2]\cup\{n+3\}$ and fails for $q>n+1+\beta_*$.
In $\mathbb R^2$, the valid range improves to $[1,5]$.
}
%Here $2<\beta_*<3$.
% The gray segment $(-1,0)$
% is shown only to indicate the extension of the chord integral
% to negative indices; no Brunn-Minkowski conclusion is asserted there.}
\label{fig-chord-ranges}
\end{figure}
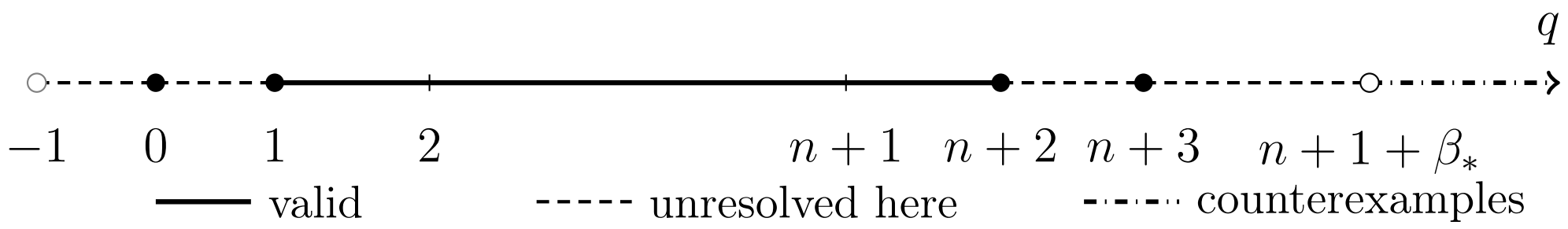

The results above establish the chord Brunn-Minkowski inequality on $[1,n+2]\cup\{n+3\}$, and in the planar case, for every $1\leq q\leq5$.
For $0<q<1$, we prove the sharp inequality for pairs of simplices and a factor inequality for arbitrary convex bodies.
The inequality fails for $q>n+1+\beta_*$. These ranges are summarized in Figure~\ref{fig-chord-ranges}.

% \caption{Known parameter ranges for the chord Brunn-Minkowski
% inequality. The diagram is schematic and not to scale.
% For general convex bodies, the inequality is proved on
% $[1,n+2]\cup\{n+3\}$ and fails for $q>n+1+\beta_*$.
% In dimension two, the valid range improves to $[1,5]$.
% Here $2<\beta_*<3$.}

These results leave several questions concerning the unresolved parameter ranges, equality cases, and the associated Minkowski problems.
The fractional comparison is sharp for simplices throughout $0<q<1$, whereas for arbitrary convex bodies we obtain only the factor inequality \eqref{eq:frac-factor-BM}.

\begin{prob}\label{prob-fractional}
Does the chord Brunn-Minkowski inequality \eqref{BM-Iq} hold for all convex bodies when $0<q<1$? In particular, can the factor in \eqref{eq:frac-factor-BM} be replaced by one?
\end{prob}

The transfer principle proves the inequality for all convex bodies up to $q=n+2$, while the inertia argument gives the isolated endpoint $q=n+3$. The rectangular perturbation threshold $n+1+\beta_*$ concerns a particular second variation and is not known to be the sharp boundary for the full problem.

\begin{prob}\label{prob-upper-range}
What is the exact validity range of \eqref{BM-Iq} in
\[
 n+2<q<n+3
 \quad\text{and}\quad
 n+3<q\leq n+1+\beta_*?
\]
In particular, does the inequality fail for every $q>n+3$?
\end{prob}

For affine-section integrals, the convex-lifting argument gives the range $1\leq q\leq1+1/m$, while the affine Blaschke-Petkantschin formula gives the additional case $q=2$.

\begin{prob}\label{prob-affine-range}
For $1<m<n$, determine the validity of \eqref{BM-Imq} outside $[1,1+1/m]\cup\{0,2\}$.
\end{prob}

Finally, we return to the uniqueness problem for chord measures raised in \cite{LXYZ,HYZ}. The results above settle ordinary chord-measure uniqueness for $q\in[1,n+2]\cup\{n+3\}$, and in the planar case for $1\leq q\leq5$. 
%They also settle the corresponding $L_p$ uniqueness statement for $p>1$ in these ranges.
The failure of the Brunn-Minkowski inequality at high indices does not by itself imply nonuniqueness of chord measures.

\begin{prob}\label{prob-chord-uniqueness}
For which remaining values of $q>0$ does
\[
 F_q(K,\cdot)=F_q(L,\cdot)
\]
imply that $K$ and $L$ differ by a translation? In particular, does ordinary chord-measure uniqueness hold for $0<q<1$ and $n+2<q<n+3$?
% \[
%  0<q<1,\qquad
%  2<q<n+1,\qquad
%  n+2<q<n+3,
% \]
% and for the unresolved high-index range below $n+1+\beta_*$?
\end{prob}

\end{document}